\documentclass[12pt]{amsart}
\usepackage[english]{babel}
\usepackage{amsmath}
\usepackage{amsthm}
\usepackage{amssymb}
\usepackage{mathrsfs}
\usepackage{enumerate}
\usepackage[notcite, final, notref]{showkeys}
\usepackage{amsfonts}
\usepackage{dsfont}
\usepackage{tikz}
\usepackage[siunitx]{circuitikz}
\usetikzlibrary{arrows,calc,chains,shapes,dsp}
\def\C{\mathbb C}
\def\R{{\mathbb R}}
\newtheorem{Pa}{Paper}[section]
\newtheorem{Tm}[Pa]{{\bf Theorem}}
\newtheorem{La}[Pa]{{\bf Lemma}}
\newtheorem{Ob}[Pa]{{\bf Observation}}
\newtheorem{Cy}[Pa]{{\bf Corollary}}
\newtheorem{Rk}[Pa]{{\bf Remark}}
\newtheorem{Pn}[Pa]{{\bf Proposition}}

\newtheorem{Ex}[Pa]{{\bf Example}}
\newtheorem{Dn}[Pa]{{\bf Definition}}
\usepackage{xcolor}

\title[Quantitatively Hyper-Positive Real Functions]
{Quantitatively Hyper-Positive Real Functions}

\author[D. Alpay]{Daniel Alpay}
\address{(DA)
Faculty of Mathematics, Physics, and Computation\\
Schmidt College of Science and Technology\\
Chapman University\\
One University Drive
Orange, California 92866\\
USA}
\email{alpay@chapman.edu}
\thanks{Daniel Alpay thanks the Foster G. and Mary McGraw Professorship in
Mathematical Sciences, which supported this research.}

\author[I. Lewkowicz]{Izchak Lewkowicz}
\address{(IL) School of Electrical and Computer Engineering\\
Ben-Gurion University of the Negev\\ P.O.B. 653\\ Beer-Sheva, 84105\\
Israel}
\email{izchak@bgu.ac.il}

\begin{document}
\bibliographystyle{plain}
\begin{abstract}
Hyper-positive real, matrix-valued, rational functions are
associated with absolute stability (the Lurie problem). Here,
quantitative subsets of Hyper-positive functions, related
through nested inclusions, are introduced.
Structurally, this family of functions turns out to be
matrix-convex and closed under inversion. 
\vskip 0.2cm

\noindent
A state-space characterization of these functions through
a corresponding Kalman-Yakubovich-Popov Lemma, is given.
Technically, the classical Linear Matrix Inclusions, associated
with passive
systems, are here substituted by Quadratic Matrix Inclusions.
\end{abstract}
\maketitle


\date{today}
\tableofcontents

\bibliographystyle{plain}
\section{Introduction}
\setcounter{equation}{0}

Positive real rational functions are known to be a model for passive,
continuous-time, linear, time-invariant systems, see e.g.
\cite[Theorem 2.6.1]{AnderVongpa1973},
\cite[Section 3.18]{Belev1968}, \cite[Section 6.3]{Khalil2000}, 
\cite{Will1976}.
In \cite{Lewk2020a} it was pointed that the family of \mbox{$m\times m$-valued}
positive real rational functions, forms a cone, which is
matrix-convex\begin{footnote}{a notion which is more strict than classical
convexity}\end{footnote} and also closed under inversion. 
Moreover, it was also shown that this property is maximal in the sense
that if one picks up a function, which is analytic in the right half plane,
but not positive real, by taking 
positive scaling, sums and inversion of this function along with
$F(s)\equiv I_m$ (a zero degree positive real function), one can always
obtain a function which is no longer analytic in the right-half plane.
\vskip 0.2cm

\noindent
It turned out that corresponding state-space realizations are also
intimately linked to the structure of matrix-convex cones, closed under
inversion. This can be intuitively explained by employing LMI (Linear Matrix
Inclusions\begin{footnote}{meaning ``inclusion within the set of positive
(semi)-definite matrices" which is more accurate than mere
``inequality".}\end{footnote}) formulation.
\vskip 0.2cm

\noindent
In this work we focus on a {\em dissipative} subset of positive real
functions\begin{footnote}{In fact, a proper subset of Strictly Positive
Real functions.}\end{footnote}: The~ {\em hyper}-positive real rational
functions, which are associated with absolute stability (a.k.a. the Lurie
problem). We then introduce a partial ordering within this set. It turns
out that here, the corresponding structure is of matrix-convex sets (not
cones), closed under inversion. Furthermore here, the LMI's are
substituted by Quadratic Matrix Inclusions.
\vskip 0.2cm

\noindent
The work is outlined as follows. In Sections \ref{sec:sub-unitDisk} and
\ref{Sect:HBeta}, we introduce a quantitative nested set of inclusions
into the Stein matrix inclusion and of the family of Bounded real
rational functions, respectively. Applying the Cayley transform, the
corresponding refinement of the Riccati matrix inclusion and the
Positive real rational functions, are given in Sections
\ref{sec:MatricesCis} and \ref{Sect:HPeta}, respectively. In Section
\ref{Sect:K-Y-P} a corresponding Kalman-Yakubovich-Popov Lemma is
presented.

\section[Sub-Unit disk]{Matrices whose spectrum is within
a subunit Disk}
\setcounter{equation}{0}
\label{sec:sub-unitDisk}

\noindent
In the sequel, we shall denote by $\overline{\mathbf H}_n$ ($\mathbf{H}_n$) the
set of $n\times n$ Hermitian (non-singular) matrices. Similarly
($\overline{\mathbf P}_n$) $\mathbf{P}_n$ will be set of $n\times n$ positive
(semi)-definite matrices\begin{footnote}{Recall, $\overline{X}$ denotes the
closure of a set $X$.}\end{footnote}.
\vskip 0.2cm

\noindent
For convenience, when using the Cayley transform, we follow the (artificial
yet intuitive) engineering notational distinction between $A$ (typically with
spectrum in $\C_R$) and~ $\hat{A}~$ (whose spectrum is typically in the unit
disk).
\vskip 0.2cm

\noindent
For a given $H\in\mathbf{H}_n$, let us define the set of all matrices $\hat{A}$
satisfying the Stein\begin{footnote}{In some engineering circles referred to a 
``discrete-time Lyapunov".}\end{footnote} inclusion,
\begin{equation}\label{eq:SteinCommonH}
\begin{matrix}
\mathbf{Stein}_H&=&\left\{ \hat{A}\in\C^{n\times n}~:~
\left(\begin{smallmatrix}\hat{A}\\~\\I_n\end{smallmatrix}\right)^*
\left(\begin{smallmatrix}
-H&0\\~\\0&H
\end{smallmatrix}\right)
\left(\begin{smallmatrix}\hat{A}\\~\\ I_n\end{smallmatrix}\right)
\in{\mathbf P}_n~\right\}.
\end{matrix}
\end{equation}
Properties of the set $\mathbf{Stein}_H$ are explored in
\cite{Lewk2020c}.
We now introduce a scalar parameter $~{\scriptstyle\eta}$, where
$~{\scriptstyle\eta}\in(1,~\infty]$, to the sets in Eq. \eqref{eq:SteinCommonH},
\begin{equation}\label{eq:SteinRhoCommonH}
\begin{matrix}
\mathbf{Stein}_H({\scriptstyle\eta})&=&\left\{ \hat{A}\in\C^{n\times n}~:~
\left(\begin{smallmatrix}\hat{A}\\~\\I_n\end{smallmatrix}\right)^*
\left(\begin{smallmatrix}
-(\eta+1)H&&0\\~\\0&&(\eta-1)H
\end{smallmatrix}\right)
\left(\begin{smallmatrix}\hat{A}\\~\\ I_n\end{smallmatrix}\right)
\in\mathbf{P}_n~\right\}
\end{matrix}
\end{equation}
Note that for all $H\in\mathbf{H}_n$ and all
$~{\scriptstyle\eta}\in(1,~\infty]$, the matrix
\[
\left(\begin{smallmatrix}
-(\eta+1)H&&0\\~\\0&&(\eta-1)H
\end{smallmatrix}\right)
\]
has real eigenvalues, $n$ of them are positive and $~n$ negative.
\vskip 0.2cm

\noindent
Eq. \eqref{eq:SteinRhoCommonH} offers a more detailed examination of
the set in Eq. \eqref{eq:SteinCommonH}, in the following sense,
\[
\infty>\eta>{\eta}_1>1\quad\Longrightarrow\quad
\mathbf{Stein}_H({\scriptstyle\eta}_1)\subset\mathbf{Stein}_H({\scriptstyle\eta})
\subset\mathbf{Stein}_H~,
\]
and
\[
\lim\limits_{\eta~\longrightarrow~\infty}\mathbf{Stein}_H({\scriptstyle\eta})
=\mathbf{Stein}_H~.
\]
Here are three basic properties of this set.

\begin{Pn}\label{Pn:ClosureOfS_H}
For parameters $~H\in\mathbf{H}_n$ and \mbox{${\scriptstyle\eta}\in(1,~\infty]$},
the matricial set $\mathbf{Stein}_H({\scriptstyle\eta})$ in Eq.
\eqref{eq:SteinRhoCommonH}, satisfies the following,

\begin{itemize}
\item[(i)~~]{}
The set $\mathbf{Stein}_H({\scriptstyle\eta})$ is open\begin{footnote}{In
the sense that $\overline{\mathbf{P}}_n$ is the closure of
$\mathbf{P}_n$.}\end{footnote} and in particular contains all elements
of the form $~{\scriptstyle\lambda}I_n$, where $~{\scriptstyle\lambda}\in\C$
and $~{\scriptstyle\sqrt{\frac{\eta-1}{\eta+1}}}>|{\scriptstyle\lambda}|$.
\vskip 0.2cm

\item[(ii)~]{}
This set is closed under multiplication by 
complex scalars of bounded modulus, i.e.
\[
\hat{A}\in\mathbf{Stein}_H({\scriptstyle\eta})\quad
\left.
\begin{smallmatrix}
{\scriptstyle\lambda}\in\C\\~\\
{\scriptstyle\sqrt{\frac{\eta-1}{\eta+1}}}>
|{\scriptstyle\lambda}|
\end{smallmatrix}
\right\}
\quad\Longrightarrow\quad
{\scriptstyle\lambda}\hat{A}\in\mathbf{Stein}_H({\scriptstyle\eta}).
\]
\item[(iii)]{}
This set is matrix-product-contractive, i.e.
\[
\hat{A}_a, \hat{A}_b\in\mathbf{Stein}_H({\scriptstyle\eta})
\quad\Longrightarrow\quad
(\hat{A}_a\hat{A}_b)\in\mathbf{Stein}_H({\scriptstyle\eta}_1)
\quad\quad{\rm where}\quad\quad{\scriptstyle{\eta}_1}
={\scriptstyle\frac{1}{2}}
({\scriptstyle\eta}+{\scriptstyle\frac{1}{\eta}}).
\]
\end{itemize}
\end{Pn}

\noindent
The above properties are easy to verify.
\vskip 0.2cm

\noindent
In
\cite[Theorem 3.5]{Ando2004} T. Ando showed that that items (i), (ii), (iii)
in Proposition \ref{Pn:ClosureOfS_H} are three out of seven
items ~{\em characterizing}~ any set of the form $\mathbf{Stein}_H$ (without
specifying $H$). The details are technically demanding and thus omitted.
\vskip 0.2cm

\noindent
In the rest of this section we shall confine the discussion to
$H\in\mathbf{P}_n$. Then, by multiplying Eq. \eqref{eq:SteinRhoCommonH},
by $H^{-\frac{1}{2}}$ from both sides, one can equivalently write
\begin{equation}\label{eq:AlternativeSteinRhoCommonH}
\begin{matrix}
\mathbf{Stein}_H({\scriptstyle\eta})&=&\{ \hat{A}\in\C^{n\times n}~:~
{\scriptstyle\sqrt{\frac{\eta-1}{\eta+1}}}
>\| H^{\frac{1}{2}}\hat{A}H^{-\frac{1}{2}}\|_2~\},
\end{matrix}
\end{equation}
where $\|~\|_2$ denotes the spectral (a.k.a. Euclidean) norm.
As a motivation, one has the following model of a
stability robustness problem.
\vskip 0.2cm

\noindent
{\bf Difference inclusion stability problem}\\
Let $x(\cdot)$ be a real vector-valued sequence satisfying,
\begin{equation}\label{eq:DifferenceIncl}
x(k+1)=A\left(k,x(k)\right)x(k)
\quad
\quad
\quad
k=0,~1,~2,~\ldots
\end{equation}
where the actual sequence $\{A\left(0,x(0)\right), ~
A\left(1,x(1)\right), ~A\left(2, x(2)\right),~\ldots~\}$ can be~ 
{\em arbitrary}.

\begin{Pn}
If in the above difference inclusion,
\[
A(\cdot , \cdot)\in\mathbf{Stein}_{I_n}({\scriptstyle\eta}),
\]
then,
\[
\|x(0)\|_2\cdot{\scriptstyle\left(\frac
{\eta-1}{\eta+1}\right)}^{\frac{k}{2}}
\geq\|x(k)\|_2\quad\quad\quad\forall~x(0)\quad
k=0,~1,~2,~\ldots,
\]
\end{Pn}
\vskip 0.2cm

\noindent
The restriction in Eq. \eqref{eq:AlternativeSteinRhoCommonH} of the set
$\mathbf{Stein}_H({\scriptstyle\eta})$ to $H\in\mathbf{P}$, enables one to
write Eq. \eqref{eq:AlternativeSteinRhoCommonH} and thus introduce an
additional property. To describe it, we need the following.

\begin{Dn}\label{Dn:MatrixConvex}
A set $\mathbf{A}$ of $n\times n$ matrices is called matrix-convex
if for all natural $k$ and all $n\times nk$ isometry
$\Upsilon=\left(\begin{smallmatrix}
{\upsilon}_1
\\
\vdots
\\~\\
{\upsilon}_k
\end{smallmatrix}\right)
$
i.e.
$~I_n={\Upsilon}^*{\Upsilon}=\sum\limits_{j=1}^k{\upsilon}_j^*{\upsilon}_j$
one has that whenever the matrices $~A_1,~\ldots~,~A_k~$ are in $\mathbf{A}$
also
\[
{\Upsilon}^*
\left(\begin{matrix}
A_1&~&~\\
&\ddots&~\\
~&~&A_k
\end{matrix}\right)
{\Upsilon}
=
\sum\limits_{j=1}^k{\upsilon}_j^*A_j{\upsilon}_j
\]
belongs to $\mathbf{A}$.
\qed
\end{Dn}
\vskip 0.2cm

\noindent
Background to matrix-convex sets in our context, is given in 
\cite[Sections 2,3]{Lewk2020a}, for broader scope literature, see
references therein. In particular, from part (I) of
\cite[Observation 2.2]{Lewk2020a} one can draw the following
conclusion with respect to the set
$\mathbf{Stein}_H({\scriptstyle\eta})$ introduced in
Eq. \eqref{eq:SteinRhoCommonH}.

\begin{Pn}\label{Stein(Eta)}
For parameters $~H\in\mathbf{P}_n~$ and a scalar
\mbox{${\scriptstyle\eta}\in(1, \infty]$} the matricial set
$\mathbf{Stein}_H({\scriptstyle\eta})$ in Eq.
\eqref{eq:AlternativeSteinRhoCommonH} is matrix-convex.
\end{Pn}
\vskip 0.2cm

\noindent
\section{Quantitatively Hyper-Bonded real rational functions}
\label{Sect:HBeta}
\setcounter{equation}{0}

\noindent
We shall use the following notation for disks in the complex plane 
\[
\mathbb{D}({\scriptstyle{\rm Center}} ,~{\scriptstyle{\rm Radius}})
=\{ {\scriptstyle\lambda}\in\C~:~{\scriptstyle{\rm Radius}}>
|{\scriptstyle\lambda}-{\scriptstyle{\rm Center}}|~\}
\quad\quad\quad\quad\begin{smallmatrix}
{\rm C}_{\rm enter}\in\C\\~\\
~{\rm R}_{\rm adius}>0.
\end{smallmatrix}
\]
In scalar terminology, a rational function is called ~{\em Bounded Real}~ if
it analytically maps $~\C_R~$ to the closed unit disk
${\mathbb D}(\underbrace{0}_{\rm Center}, \underbrace{1}_{\rm Radius})$.
Formally, a \mbox{$m\times m$-valued}
rational function $~G(s)~$ is called ~{\em Bounded Real} if
\begin{equation}\label{eq:BoundedRealDef}
\left(I_m-
\left(G(s)\right)^*G(s)\right)
\in\overline{\mathbf P}_m
\quad\quad\quad\forall s\in\C_R~,
\end{equation}
see e.g. \cite[Section 2.6]{AnderVongpa1973} and further analysis in
\cite[Section 3]{Lewk2020c}.
\vskip 0.2cm

\noindent
We next ``zoom into" this family. To simplify the exposition, we start with
scalar rational functions $g(s)$ mapping
$\C_R$ to a closed sub-unit disk of the form.
\[
{\mathbb D}(\underbrace{0}_{\rm Center},
\underbrace{\scriptstyle\sqrt{\frac{\eta-1}{\eta+1}}}_{\rm Radius})
\quad\quad\quad{\scriptstyle\eta}\in(1,~\infty],
\]
This set of $\eta$-Hyper Bounded Real functions, will be denoted by
$\mathcal{HB}_{\eta}~$. For
\mbox{$m\times m$-valued} rational functions, we shall
describe this family as
subsets of functions $G(s)$ in Eq. \eqref{eq:BoundedRealDef},
\begin{equation}\label{eq:EtaBoundedRealFunction}
\left(
I_m-{\scriptstyle\frac{\eta+1}{\eta-1}}\left(G(s)\right)^*G(s)\right)
\in\overline{\mathbf P}_m
\quad\quad\quad\forall s\in\C_R~.
\end{equation}
This condition can be equivalently written as
\begin{equation}\label{eq:NormBound}
\sqrt{\scriptstyle\frac{\eta-1}{\eta+1}}\leq\sup\limits_{s\in\C_R}\|G(s)\|_2~.
\end{equation}
In the sequel we shall illustrate the contractive nature of the
sets of Hyper-Bounded (and subsequently Hyper-Positive) in several
ways, here are scalar samples of $\mathcal{HB}_{\eta}$ functions.

\begin{Ex}\label{Ex:scalarHBeta}
{\rm
Consider a degree one $\mathcal{HB}_{\eta}$ rational function $g(s)$ of the
form
\begin{equation}\label{eq:DefScalarHyperBR}
g(s)={\scriptstyle\sqrt{\frac{\eta-1}{\eta+1}}}\cdot\frac{s-a}{s+a}
\quad\quad\begin{smallmatrix}a>0\\~\\
\eta\in(1,~\infty].\end{smallmatrix}
\end{equation}
It maps $\C_R$ ~{\em onto}~ the closed sub-unit disk $\mathbb{D}(
\underbrace{0}_{\rm Center}, \underbrace{\scriptstyle\frac{\sqrt{\eta-1}}
{\sqrt{\eta+1}}}_{\rm Radius})$, irrespective of the parameter $~a$.  In particular,
\mbox{$~g(0)=-{\scriptstyle\sqrt{\frac{\eta-1}{\eta+1}}}$},
\mbox{$~g(ia)=-i{\scriptstyle\sqrt{\frac{\eta-1}{\eta+1}}}$},
\mbox{$~g(\infty)={\scriptstyle\sqrt{\frac{\eta-1}{\eta+1}}}$}. See the Nyquist plot
on the left-hand side of Figure \ref{Figure:DiskCISs2}.
\qed
}
\end{Ex}
\vskip 0.2cm

\noindent
We now address ourselves to state space realizations of the above rational
function $G(s)$. We start with the classical
Bounded Real Lemma, see e.g. \cite[Section 7.2]{AnderVongpa1973},
\cite[Subsection 2.7.3]{BGFB1994}.

\begin{La}\label{La:ClassicalBRlemma}
Let $G(s)$ be a $m\times m$-valued rational function, with no pole
at infinity, of McMillan degree $n$. Let
\[
G(s)=\hat{C}(sI_n-\hat{A})^{-1}\hat{B}+\hat{D},
\]
be a corresponding minimal realization.
\vskip 0.2cm

\noindent
The function $G(s)$ is bounded real, if and only if,
there exist a matrix $H$ satisfying
\[
H\in\mathbf{P}_n\quad\quad\quad
\left(
\left(\begin{smallmatrix}-H\hat{A}-\hat{A}^*H
&~&-H\hat{B}\\ -\hat{B}^*H&~&I_m\end{smallmatrix}\right)
-
\left(
\begin{smallmatrix}
\hat{C}^*\\ \hat{D}^*
\end{smallmatrix}
\right)
\left(
\begin{smallmatrix}
\hat{C}&\hat{D}
\end{smallmatrix}\right)
\right)\in\overline{\mathbf P}_{n+m}~.
\]
\end{La}
\vskip 0.2cm

\noindent
One can now introduce a quantitative refinement of this classical result.

\begin{La}\label{La:QuantitativeBRlemma}
Let $G(s)$ be a \mbox{$m\times m$-valued} rational function, with no pole
at infinity, of McMillan degree $~n$. Let
\[
G(s)=\hat{C}(sI_n-\hat{A})^{-1}\hat{B}+\hat{D},
\]
be a corresponding minimal realization.
\vskip 0.2cm

\noindent
We shall say that $m\times m$-valued rational function, belongs to
$\mathcal{HB}_{\eta}$ for some \mbox{${\scriptstyle\eta}\in(1,~\infty]$},
i.e.  $G(s)$ is quantitatively bounded real in the sense of Eq.
\eqref{eq:EtaBoundedRealFunction}, if and only if, there exist a matrix
$H$ satisfying
\[
H\in\mathbf{P}_n\quad\quad\quad\left(
\left(\begin{smallmatrix}-H\hat{A}-\hat{A}^*H
&~&-H\hat{B}\\ -\hat{B}^*H&~&I_m\end{smallmatrix}\right)
-
{\scriptstyle\frac{\eta+1}{\eta-1}}
\left(
\begin{smallmatrix}
\hat{C}^*\\ \hat{D}^*
\end{smallmatrix}
\right)
\left(
\begin{smallmatrix}
\hat{C}&\hat{D}
\end{smallmatrix}\right)
\right)\in\overline{\mathbf P}_{n+m}~.
\]
Moreover, the right-hand side is singular.
\end{La}
\vskip 0.2cm

\noindent
The classical case (as in Lemma \ref{La:ClassicalBRlemma})
is obtained upon substituting $~{\scriptstyle\eta}=\infty$.
\vskip 0.2cm

\noindent
In the sequel we shall find it convenient to re-write the condition in Lemma
\ref{La:QuantitativeBRlemma} in the following quadratic form.
\begin{equation}\label{eq:QuantitativeBRlemmaAlternative}
H\in\mathbf{P}_n\quad\quad\quad
\left(\begin{smallmatrix}
\hat{A}&\hat{B}\\
\hat{C}&\hat{D}\\
I_n&0\\
0&I_m
\end{smallmatrix}\right)^*
\underbrace{
\left(\begin{smallmatrix}
~~0&0&-H&~0\\
~~0&{\scriptstyle\frac{1+\eta}{1-\eta}}I_m&~~0&~0\\
-H&0&~~0&~0\\
~~0&0&~~0&~I_m
\end{smallmatrix}\right)
}_W
\left(\begin{smallmatrix}
\hat{A}&\hat{B}\\
\hat{C}&\hat{D}\\
I_n&0\\
0&I_m
\end{smallmatrix}\right)\in\overline{\mathbf P}_{n+m}~.
\end{equation}
Note that the Hermitian matrix $W$,
has exactly $n+m$ positive eigenvalues and $n+m$ negative eigenvalues.
\vskip 0.2cm

\noindent
We next resort to the classical Cayley transform.

\begin{Dn}\label{Dn:MatrixCayleyTransform}
{\rm
We denote by $\mathcal{C}(M)$ the Cayley transform of a matrix $M\in\C^{n\times n}$,
\begin{center}
$\begin{matrix}
\hat{M}&=&\mathcal{C}\left(M\right)&:=&\left(I_n-M\right)\left(I_n+M\right)^{-1}
=-I_n+2\left(I_n+M\right)^{-1}
&~&~&~&-1\not\in{\rm spect}(M).
\end{matrix}$
\end{center}
\qed
}
\end{Dn}

\noindent
As already mentioned, for convenience we follow the
(artificial yet intuitive)
engineering notational distinction between $M$ (typically with
spectrum in $\C_R$) and~ $\hat{M}=\mathcal{C}(M)$ ~(whose spectrum is
typically in the unit disk).
\vskip 0.2cm

\noindent
Recall that the Cayley transform is involutive in the sense that,
whenever well defined,
\[
\mathcal{C}\left(\mathcal{C}\left(M\right)\right)=M.
\]
\vskip 0.2cm

\begin{center}
\begin{figure}
\begin{minipage}[b]{0.47\linewidth}
\begin{tikzpicture}[scale=3,cap=round]
    \tikzstyle{axes}=[]
     \tikzstyle{important line}=[very thick]
     \tikzstyle{information text}=[rounded corners,fill=red!10,inner sep=1ex]
     \begin{scope}[style=axes]
   \
       \draw[->] (-0.7,0) -- (0.7,0) node[right] {Real};
       \draw[->] (0,-0.7) -- (0,0.7) node[above] {Imaginary};

       \foreach \x/\xtext in
{- 0.5/-\frac{\sqrt{\eta-1}}{\sqrt{\eta+1}},
0.5/\frac{\sqrt{\eta-1}}{\sqrt{\eta+1}}
}
\draw[xshift=\x cm] (0pt,1pt) -- (0pt,-1pt) node[below,fill=white]
    {$\xtext$}; %If we remove the node then we won't have the numbering below the grid

       \foreach \y/\ytext in
{- 0.5/-\frac{\sqrt{\eta-1}}{\sqrt{\eta+1}}, 0.5/\frac{\sqrt{\eta-1}}{\sqrt{\eta+1}}
}
        \draw[yshift=\y cm] (1pt,0pt) -- (-1pt,0pt) node[left,fill=white]
     {$\ytext$};%If we remove the node then we won't have the numbering below the grid
     \end{scope}
    \draw[arrows=->,style=important line, red] (0,0) circle (0.5);
\end{tikzpicture}
\label{Figure:SubUnitDisks2}
\mbox{Nyquist plot of 
$g(s)$
in Eq. \eqref{eq:DefScalarHyperBR}}
\end{minipage}
\quad\quad
\begin{minipage}[b]{0.47\linewidth}
\begin{tikzpicture}[scale=1.2,cap=round]
     \tikzstyle{axes}=[]
     \tikzstyle{important line}=[very thick]
     \tikzstyle{information text}=[rounded corners,fill=red!10,inner sep=1ex]
     \begin{scope}[style=axes]
   \
       \draw[->] (-0.5,0) -- (3.6,0) node[right] {Real};
       \draw[->] (0,-1.6) -- (0,1.7) node[above] {Imaginary};

       \foreach \x/\xtext in {
0.333/{\scriptstyle\eta-\sqrt{{\eta}^2-1}}, 1.67/{\scriptstyle\eta},
3/{\scriptstyle\eta+\sqrt{{\eta}^2-1}}}
 \draw[xshift=\x cm] (0pt,1pt) -- (0pt,-1pt) node[below,fill=white]
     {$\xtext$}; %If we remove the node then we won't have the numbering below the grid

    \foreach \y/\ytext in {-1.33/{\scriptstyle-\sqrt{{\eta}^2-1}},
1.33/{\scriptstyle\sqrt{{\eta}^2-1}}}
        \draw[yshift=\y cm] (1pt,0pt) -- (-1pt,0pt) node[left,fill=white]
       {$\ytext$};%If we remove the node then we won't have the numbering below the g
    \end{scope}
    \draw[arrows=->,style=important line, red] (5/3,0) circle (4/3);
  \end{tikzpicture}
\mbox{
Nyquist plot of 
$f(s)$ in Eq. \eqref{eq:ExampleStateSpaceDeg1diskPagain}
}
\end{minipage}
\caption{$f(s)=\mathcal{C}\left(g(s)\right)$
}
\label{Figure:DiskCISs2}
\end{figure}
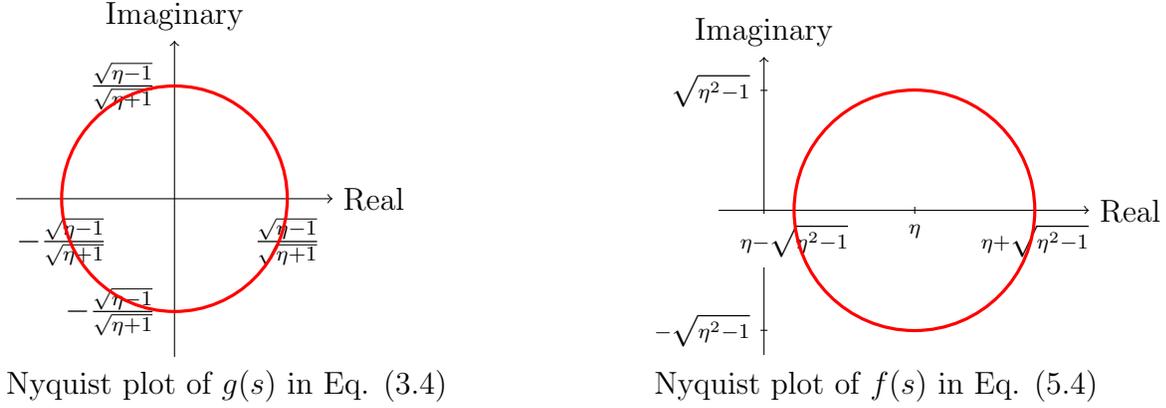
\end{center}
\vskip 0.2cm

\noindent
The contractive nature of $\mathcal{HP}_{\eta}$ functions will be further
discussed in Observation \ref{Ob:Composition} and through state-space
realization presentation, in items (2), (3) of
Example \ref{Ex:KypQuantitativeHyperStable}.
\vskip 0.2cm

\noindent
Disks of the form of the right-hand side of Figure \ref{Figure:DiskCISs2}
are studied in details in \cite{Lewk2020b}.

\section{Quantitatively contractive Lyapunov inclusions}
\label{sec:MatricesCis}
\setcounter{equation}{0}

\noindent
Applying the Cayley transform to the Stein inclusion in Eq.
\eqref{eq:SteinRhoCommonH}, yields the Riccati inclusion
in Eq. \eqref{eq:DefLh(eta)} below.

\begin{Pn}\label{Pn:CayleyLyapStein}
For arbitrary $H\in\mathbf{H}_n$ and ${\scriptstyle\eta}\in(1,~\infty]$
one has that
\begin{equation}\label{eq:CayleyStein}
\mathcal{C}\left(\mathbf{Stein}_H({\scriptstyle\eta})\right)
=\mathbf{L}_H({\scriptstyle\eta}),
\end{equation}
where,
\begin{equation}\label{eq:DefLh(eta)}
\mathbf{L}_H({\scriptstyle\eta})=
\left\{ A\in\C^{n\times n}~:~
\left(\begin{smallmatrix}A\\~\\I_n\end{smallmatrix}\right)^*
\left(\begin{smallmatrix}
-\frac{1}{\eta}H&H\\~\\H&-\frac{1}{\eta}H
\end{smallmatrix}\right)
\left(\begin{smallmatrix}A\\~\\ I_n\end{smallmatrix}\right)
=Q\in{\mathbf P}_n~\right\}.
\end{equation}
\end{Pn}
\vskip 0.2cm

\begin{Rk}\label{Rk:SteinLyapRhoCommonH}
{\rm
Note that in the Stein equation
\eqref{eq:SteinRhoCommonH},
\eqref{eq:CayleyStein} and in the Riccati equation
\eqref{eq:DefLh(eta)}, the parameters $H$ and ${\scriptstyle\eta}$, are
indeed the same.
}
\qed
\end{Rk}

\noindent
In Eq. \eqref{eq:DefLh(eta)} in particular we denote
$\lim\limits_{\eta~\longrightarrow~\infty}\mathbf{L}_H({\scriptstyle\eta})
=\mathbf{L}_H$ namely a Lyapunov inclusion,
\[
\mathbf{L}_H=
\left\{ A\in\C^{n\times n}~:~
\left(\begin{smallmatrix}A\\~\\I_n\end{smallmatrix}\right)^*
\left(\begin{smallmatrix}
0&H\\~\\H&0
\end{smallmatrix}\right)
\left(\begin{smallmatrix}A\\~\\ I_n\end{smallmatrix}\right)
=HA+A^*H=Q\in{\mathbf P}_n~\right\}.
\]
For a prescribed $~H\in\mathbf{H}_n$ one has that,
\[
\infty>{\scriptstyle\eta}>{\scriptstyle{\eta}_1}
>1
\quad\Longrightarrow\quad
\mathbf{L}_H({\scriptstyle{\eta}_1})
\subset
\mathbf{L}_H({\scriptstyle\eta})
\subset
\mathbf{L}_H~,
\]
and each inclusion is strict.
\vskip 0.2cm

\noindent
The sets $\overline{\mathbf L}_H$, $\mathbf{L}_H$ (where
${\scriptstyle\eta}=\infty$) were introduced in \cite{CohenLew1997a},
where it was shown that they are convex cones and closed under inversion.
In \cite{Lewk2020a} it was shown that for $H=I_n$, these sets are in fact
{\em matrix-convex}.  Here, we colloquially ``zoom-in" into this structure
to present the main result of this section.

\begin{Tm}\label{Tm:ClosureOfL_H(r)}
Let $~H\in\mathbf{P}_n~$ and a scalar $~{\scriptstyle\eta}\in(1,~\infty]$,
be prescribed.\\
(I) Then, the open set $~\mathbf{L}_H({\scriptstyle\eta})$ in Eq.
\eqref{eq:DefLh(eta)} is convex and closed under inversion.
\vskip 0.2cm

(II) For  $~H\in\mathbf{P}_n~$ and ${\scriptstyle\eta}=\infty$ the set
$~\mathbf{L}_H$ is a maximal convex set of matrices whose spectrum is in
$\C_R~$.
\vskip 0.2cm

(III) For $~H=I_n$, the set $~{\mathbf L}_{I_n}({\scriptstyle\eta})$, is 
in addition matrix-convex.
\end{Tm}

\noindent
{\bf Proof :}\quad 
(I)\quad To show convexity, for a pair of $n\times n$ matrices
$A_o$, $A_1$ and a scalar parameter
\mbox{${\scriptstyle\theta}\in[0,~1]$},
consider the following identity, 
\[
\begin{matrix}
{\scriptstyle\theta}
\underbrace{
\left(\begin{smallmatrix}A_1\\~\\I_n\end{smallmatrix}\right)^*
\left(\begin{smallmatrix}
-\frac{1}{\eta}H&H\\~\\H&-\frac{1}{\eta}H
\end{smallmatrix}\right)
\left(\begin{smallmatrix}A_1\\~\\ I_n\end{smallmatrix}\right)
}_{Q_1}
+
\left({\scriptstyle 1}-{\scriptstyle\theta}\right)
\underbrace{
\left(\begin{smallmatrix}A_o\\~\\I_n\end{smallmatrix}\right)^*
\left(\begin{smallmatrix}
-\frac{1}{\eta}H&H\\~\\H&-\frac{1}{\eta}H
\end{smallmatrix}\right)
\left(\begin{smallmatrix}A_o\\~\\ I_n\end{smallmatrix}\right)
}_{Q_o}
+
{\scriptstyle\frac{\theta(1-\theta)}{\eta}}
\left({\scriptstyle A_o-A_1}\right)^*
{\scriptstyle H}
\left({\scriptstyle A_o-A_1}\right)
\\~\\
=
\left(\begin{smallmatrix}\theta{A}_1+(1-\theta)A_o\\~\\I_n\end{smallmatrix}\right)^*
\left(\begin{smallmatrix}
-\frac{1}{\eta}H&H\\~\\H&-\frac{1}{\eta}H
\end{smallmatrix}\right)
\left(\begin{smallmatrix}\theta{A}_1+(1-\theta)A_o\\~\\ I_n\end{smallmatrix}\right).
\end{matrix}
\]
The condition $~A_o, A_1\in\mathbf{L}_H({\scriptstyle\eta})$ is equivalent to
having $Q_o, Q_1\in\mathbf{P}_n$. Since by assumption  $H\in\mathbf{P}_n$, it
implies that ~$\left({\scriptstyle A_o-A_1}\right)^*{\scriptstyle H}\left(
{\scriptstyle A_o-A_1}\right)\in\overline{\mathbf P}_n$ and thus one has
that for all \mbox{${\scriptstyle\theta}\in[0,~1]$},
also $\left({\scriptstyle\theta}{A}_1+(1-{\scriptstyle\theta})A_o\right)$
belongs to $~\mathbf{L}_H({\scriptstyle\eta})$,
so this part of the claim is established.
\vskip 0.2cm

\noindent
We shall show invertibility in two ways:
\vskip 0.2cm

\noindent
(a)~~Directly:~Using Eq. \eqref{eq:DefLh(eta)} note that whenever
$A\in\mathbf{L}_H({\scriptstyle\eta})$ then,
\[
\left(\begin{smallmatrix}A^{-1}\\~\\I_n\end{smallmatrix}\right)^*
\left(\begin{smallmatrix}
-\frac{1}{\eta}H&H\\~\\H&-\frac{1}{\eta}H
\end{smallmatrix}\right)
\left(\begin{smallmatrix}A^{-1}\\~\\ I_n\end{smallmatrix}\right)
=
\left(\begin{smallmatrix}I_n\\~\\A^{-1}\end{smallmatrix}\right)^*
\left(\begin{smallmatrix}
-\frac{1}{\eta}H&H\\~\\H&-\frac{1}{\eta}H
\end{smallmatrix}\right)
\left(\begin{smallmatrix}I_n\\~\\A^{-1}\end{smallmatrix}\right)
=
{\scriptstyle\left(A^{-1}\right)^*Q\left(A^{-1}\right)
}
\in\mathbf{P}_n~.
\]
(b)~~Invertibility can be indirectly deduced by using the relation
$\mathcal{C}\left(A^{-1}\right)=-\mathcal{C}\left(A\right)$ 
which is immediate from Definition \ref{Dn:MatrixCayleyTransform},
along with Eqs. \eqref{eq:SteinRhoCommonH}
and \eqref{eq:AlternativeSteinRhoCommonH}.
\vskip 0.2cm

\noindent
(II)~~ See e.g. \cite[Lemma 3.5]{CohenLew1997a}.
\vskip 0.2cm

(III)~~Matrix convexity for $H=I_n~$.
\vskip 0.2cm

\noindent
First note that one can always complete a given $kn\times n$ isometry
$\Upsilon$ to a $kn\times kn$ unitary matrix $U$, i.e. to
construct a $kn\times(k-1)n~$ matrix $\tilde{U}$ so that
\begin{equation}\label{eq:Completion}
I_{kn}=
\underbrace{\left(\begin{matrix}\Upsilon&\tilde{U}\end{matrix}\right)}_{U}
\underbrace{\left(\begin{matrix}\Upsilon&\tilde{U}\end{matrix}\right)^*}_{U^*}
=
{\Upsilon}{\Upsilon}^*+\tilde{U}\tilde{U}^*
\end{equation}
For a prescribed ${\scriptstyle\eta}\in(1,~\infty]$ let now
$A_1,~\ldots~,~A_k$ be in
$\mathbf{L}_{I_n}({\scriptstyle\eta})$, then by Eq.\eqref{eq:DefLh(eta)},
\begin{equation}\label{eq:Auxilliary}
-{\scriptstyle\frac{1}{\eta}}A_j^*A_j+A_j^*+A_j
-{\scriptstyle\frac{1}{\eta}}I_n=
\Delta_j
\quad\quad
\begin{smallmatrix}\Delta_j\in\mathbf{P}_n\\~\\
j=1,~\ldots~,~k.
\end{smallmatrix}
\end{equation}
We need to show that each matrix of the form
\[
\sum\limits_{j=1}^k{\upsilon}_j^*{\scriptstyle A}_j{\upsilon}_j
=
{\scriptstyle\Upsilon}^*
\left(\begin{smallmatrix}
A_1&~&~\\
~&\ddots&~\\
~&~&A_k
\end{smallmatrix}\right)
{\scriptstyle\Upsilon},
\]
satisfies the condition in Eq. \eqref{eq:DefLh(eta)}.
\vskip 0.2cm

\noindent
Indeed,
\[
\begin{matrix}
~&
-{\scriptstyle\frac{1}{\eta}}
\left(\sum\limits_{j=1}^k{\upsilon}_j^*{\scriptstyle A_j}{\upsilon}_j\right)^*
\left(\sum\limits_{j=1}^k{\upsilon}_j^*{\scriptstyle A_j}{\upsilon}_j\right)
+
\left(\sum\limits_{j=1}^k{\upsilon}_j^*{\scriptstyle A_j}{\upsilon}_j\right)^*
+
\left(\sum\limits_{j=1}^k{\upsilon}_j^*{\scriptstyle A_j}{\upsilon}_j\right)
-{\scriptstyle\frac{1}{\eta}I_n}
\\
=&
-{\scriptstyle\frac{1}{\eta}}
\left(\sum\limits_{j=1}^k{\upsilon}_j^*{\scriptstyle A_j}{\upsilon}_j\right)^*
\left(\sum\limits_{j=1}^k{\upsilon}_j^*{\scriptstyle A_j}{\upsilon}_j\right)
+
\sum\limits_{j=1}^k{\upsilon}_j^*\underbrace{\left(
{\scriptstyle A}_j+{\scriptstyle A}_j^*-{\scriptstyle\frac{1}{\eta}I_n}
\right)}_{\frac{1}{\eta}A_j^*A_j+{\Delta}_j~{\rm by~Eq.} 
\eqref{eq:Auxilliary}}{\upsilon}_j
\\
=&
-{\scriptstyle\frac{1}{\eta}}
\left(\sum\limits_{j=1}^k{\upsilon}_j^*{\scriptstyle A_j}{\upsilon}_j\right)^*
\left(\sum\limits_{j=1}^k{\upsilon}_j^*{\scriptstyle A_j}{\upsilon}_j\right)
+
\sum\limits_{j=1}^k{\upsilon}_j^*
\left(
{\scriptstyle\frac{1}{\eta}A_j^*A_j}+{\scriptstyle\Delta_j}
\right){\upsilon}_j
\\
=&
-{\scriptstyle\frac{1}{\eta}}
\left(\begin{smallmatrix}
A_1{\upsilon}_1
\\
\vdots
\\
A_k{\upsilon}_k
\end{smallmatrix}
\right)^*
\Upsilon{\Upsilon}^*
\left(\begin{smallmatrix}
A_1{\upsilon}_1
\\
\vdots
\\
A_k{\upsilon}_k
\end{smallmatrix}
\right)
+
{\scriptstyle\frac{1}{\eta}}
\left(\begin{smallmatrix}
A_1{\upsilon}_1
\\
\vdots
\\
A_k{\upsilon}_k
\end{smallmatrix}
\right)^*
\left(\begin{smallmatrix}
A_1{\upsilon}_1
\\
\vdots
\\
A_k{\upsilon}_k
\end{smallmatrix}
\right)
+{\Upsilon}^*
\left(\begin{smallmatrix}
{\Delta}_1&~&~\\
~&\ddots&~\\
~&~&{\Delta}_k
\end{smallmatrix}\right)
{\Upsilon}
\\
=&
{\scriptstyle\frac{1}{\eta}}
\left(\begin{smallmatrix}
A_1{\upsilon}_1
\\
\vdots
\\
A_k{\upsilon}_k
\end{smallmatrix}
\right)^*
\underbrace{\left(I_{nk}-
\Upsilon{\Upsilon}^*
\right)}_{
\tilde{U}{\tilde{U}}^*~
{\rm by~Eq.}~\eqref{eq:Completion}
}
\left(\begin{smallmatrix}
A_1{\upsilon}_1
\\
\vdots
\\
A_k{\upsilon}_k
\end{smallmatrix}
\right)
+{\Upsilon}^*
\left(\begin{smallmatrix}
{\Delta}_1&~&~\\
~&\ddots&~\\
~&~&{\Delta}_k
\end{smallmatrix}\right)\Upsilon
\\
=&
{\scriptstyle\frac{1}{\eta}}
\left(\begin{smallmatrix}
A_1{\upsilon}_1
\\
\vdots
\\
A_k{\upsilon}_k
\end{smallmatrix}
\right)^*
\tilde{U}{\tilde{U}}^*
\left(\begin{smallmatrix}
A_1{\upsilon}_1
\\
\vdots
\\
A_k{\upsilon}_k
\end{smallmatrix}
\right)
+{\Upsilon}^*
\left(\begin{smallmatrix}
{\Delta}_1&~&~\\
~&\ddots&~\\
~&~&{\Delta}_k
\end{smallmatrix}\right)\Upsilon.
\end{matrix}
\]
As we have a sum of positive semi-definite matrices, the matrix
$\sum\limits_{j=1}^k{\upsilon}_j^*{\scriptstyle A}_j{\upsilon}_j$ satisfies
Eq. \eqref{eq:DefLh(eta)} and thus belongs to
$\mathbf{L}_{I_n}({\scriptstyle\eta})$, so the claim is established.
\qed
\vskip 0.2cm

\noindent
The above structural properties will be exploited in Corollary
\ref{Cy:AlternativeHPeta} below, where we address \mbox{$m\times m$-valued}
rational functions mapping $\C_R$ to
\[
\overline{\mathbf L}_{I_m}=
\left\{A\in\C^m~:~A+A^*\in\overline{\mathbf P}_m~\right\}.
\] 

\section{Quantitatively Hyper-Positive Real Functions}
\label{Sect:HPeta}
\setcounter{equation}{0}

We start by setting up the framework for positive real rational functions
and their subsets.
\vskip 0.2cm

\begin{Dn}
{\rm
Consider the following families of $m\times m$-valued rational
functions $F(s)$.

\begin{itemize}
\item{}The family of functions analytically mapping $\C_R$ to
$\overline{\mathbf L}_{I_m}$ is called Positive Real and denoted
by $~\mathcal{P}$.
\vskip 0.2cm

\item{}The family of functions $F(s)$ so that $F(s-\epsilon)$ is
positive real, for some $\epsilon>0$, is called Strictly Positive Real
and denoted by $~\mathcal{SP}$.
\vskip 0.2cm

\item{}The family of functions analytically mapping $\C_R$ to
${\mathbf L}_{I_m}$ is called Hyper-Positive Real and denoted
by $~\mathcal{HP}$.
\end{itemize}
}
\qed
\end{Dn}

As a scalar illustration consider the following degree one $\mathcal{P}$
unction
\[
f(s)={\scriptstyle d}+{\scriptstyle\frac{b}{s+a}}
\quad\quad\quad\begin{smallmatrix}
b>0\\~\\
a, d\geq 0~~{\rm parameters}.
\end{smallmatrix}
\]
$f\in\mathcal{SP}~$ if and only if $~ab>0$, $d\geq 0$ and
$f\in\mathcal{HP}~$ if and only if $~abd>0$.
\vskip 0.2cm

\noindent
While the sets $\mathcal{P}$ and $\mathcal{SP}$ are classical, see
e.g.  \cite{BroLozaMasEge2007}, the name
$\mathcal{HP}$, is less conventional. For instance, in some research
circles, see e.g. \cite{BroLozaMasEge2007}, this family of functions is
referred to as ``both input and output strictly passive".
\vskip 0.2cm

\noindent
Here is a summary of the structure of the set $\mathcal{P}$. For details
see \cite[Theorem 4.3]{Lewk2020a}.
\vskip 0.2cm

\begin{Ob}\label{Ob:PropertiesP}
Consider the sets of $m\times m$-valued rational functions: $\mathcal{P}$,
$\mathcal{SP}$ and $\mathcal{HP}$ 

\begin{itemize}
\item[(i)~~~]{}
$\mathcal{HP}\subset\mathcal{SP}\subset\mathcal{P}$ and each inclusion
is strict.
\vskip 0.2cm

\item[(ii)~~]{}
Each of these three sets is a matrix-convex cone and closed under inversion.
\vskip 0.2cm

\item[(iii)~]{} 
In item (ii), the set $\mathcal{P}$ is maximal in the sense that if a real
\mbox{$m\times m$-valued} rational function is not positive real, by
taking positive scaling, sums and inversion of this function, along with
$F(s)\equiv I_m$ (a zero degree positive real function), one can always
obtain a function which is no longer analytic in $\C_R$.
\end{itemize}
\end{Ob}
\vskip 0.2cm

\noindent
Recall that in \cite{Brune1} Otto Brune showed that the 
driving point immittance of a lumped \mbox{$R-L-C$} electrical network
belongs to $~\mathcal{P}$ and that an arbitrary positive real rational
function can be realized as the driving
point immittance of a lumped \mbox{$R-L-C$} electrical network. For
further details see e.g. \cite{AnderVongpa1973}, \cite{Belev1968}
and \cite{Wohl1969}. For a recent comprehensive account of circuits describing
$~\mathcal{P}$ functions of degree two, see \cite{MorelliSmith2019}.
\vskip 0.2cm

\noindent
Roughly speaking, $\mathcal{HP}$ functions correspond to the driving
point immittance of lumped \mbox{$R-L-C$} electrical networks with
no branch which is puerly reactive (i.e. without a
resistive) part.
\vskip 0.2cm

\noindent
In the sequel we refine the above definition of $\mathcal{HP}$
functions by introducing a parameter ${\scriptstyle\eta}$, so that
\[
\infty>{\scriptstyle\eta}>{\scriptstyle{\eta}_1}>1
\quad\Longrightarrow\quad
\mathcal{HP}_{{\eta}_1}
\subset
\mathcal{HP}_{\eta}
\subset
\mathcal{HP},
\]
and each inclusion is strict. Here are the details.
\vskip 0.2cm

\noindent
It is well known, see e.g. \cite[Example 2.7.1]{AnderVongpa1973},
that the families of Bounded Real
and Positive Real rational functions are related through
the Cayley transform, i.e.
\[
\mathcal{C}\left(\mathcal{B}\right)=\mathcal{P}.
\]
Following this intuition, we shall similarly describe the set
$\mathcal{HP}_{\eta}$ of Quantitatively Hyper-Positive Real,
\[
\mathcal{C}\left(\mathcal{HB}_{\eta}\right)=\mathcal{HP}_{\eta}~.
\]
Proposition \ref{Pn:CayleyLyapStein} relates, through the Cayley
transform, the matricial families
$\left(\mathbf{Stein}_H({\scriptstyle\eta})\right)$
and $\mathbf{L}_H({\scriptstyle\eta})$. In an analogous way, one has
the following for matrix-valued, rational $\mathcal{HP}_{\eta}$ 
functions.

\begin{Ob}\label{Ob:CayleyQuantitativeHyperPositive1}
Let $F(s)$ be a $m\times m$-valued rational function and
let ${\scriptstyle\eta}$, ${\scriptstyle\eta}\in(1,~\infty]$,
be a parameter.\quad
The function $F(s)$ belongs to $\mathcal{HP}_{\eta}$,
if and only if
\begin{equation}\label{eq:DefHPetaAgain1}
\left(\begin{matrix}F\\ I\end{matrix}\right)^*
\left(\begin{matrix}-{\scriptstyle\frac{1}{\eta}}I_n&I_n
\\I_n&-{\scriptstyle\frac{1}{\eta}}I_n\end{matrix}\right)
\left(\begin{matrix}F\\ I\end{matrix}\right)
\in\overline{\mathbf P}_m\quad\forall s\in\C_R~.
\end{equation}
\end{Ob}

In a way similar to the discussion following Remark
\ref{Rk:SteinLyapRhoCommonH}, taking the limit
${\scriptstyle\eta}~\longrightarrow~\infty$, the quadratic
inclusion in Eq. \eqref{eq:DefHPetaAgain1} degenerates
to a linear inclusion and the set
\[
\mathcal{P}:=\{~F(s):~F+F^*\in\overline{\mathbf P}_n~~
\forall s\in\C_R~\},
\]
of Positive Real functions, is recovered. For more details on
this set see \cite{Lewk2020a}.
\vskip 0.2cm
 
\noindent
Denoting by $\rho(M)$ the spectral radius of a square matrix $M$,
one obtains the following equivalent, and technically more
convenient, characterization of this set.

\begin{Cy}\label{Cy:QmiQuantitativeHyperPositive1}
Let $F(s)$ be a $m\times m$-valued rational function and
let ${\scriptstyle\eta}$, ${\scriptstyle\eta}\in(1,~\infty]$,
be a parameter.\quad
The function $F(s)$ belongs to $\mathcal{HP}_{\eta}$,
if and only if
\begin{equation}\label{eq:DefHPetaAgain2}
{\scriptstyle\eta}\geq
\sup\limits_{s\in\C_R}\rho\left(
\left(F^*F+I_m\right)\left(F^*+F\right)^{-1}\right).
\end{equation}
\end{Cy}

\noindent
{\bf Proof :}\quad
From Eq. \eqref{eq:DefHPetaAgain1} one has the following
chain of relations
\[
\begin{matrix}
\left(
\left(F+F^*\right)
-{\scriptstyle\frac{1}{\eta}}\left(F^*F+I\right)
\right)
&\in\overline{\mathbf P}_m
&&\forall s\in\C_R
\\
\left(
{\scriptstyle\eta}\left(F^*+F\right)-
\left(I_m+F^*F\right)
\right)
&\in\overline{\mathbf P}_m
&&\forall s\in\C_R
\\
\left(
{\scriptstyle\eta}I_m-
\left(F^*+F\right)^{-\frac{1}{2}}
\left(I_m+F^*F\right)
\left(F^*+F\right)^{-\frac{1}{2}}
\right)
&\in\overline{\mathbf P}_m
&&\forall s\in\C_R
\\
{\scriptstyle\eta}
\geq
\|
\left(F^*+F\right)^{-\frac{1}{2}}
\left(I_m+F^*F\right)
\left(F^*+F\right)^{-\frac{1}{2}}
\|_2
&
&&\forall s\in\C_R~,
\end{matrix}
\]
which in turn means that
\[
\begin{matrix}
{\scriptstyle\eta}&\geq&
\sup\limits_{s\in\C_R}
\|\left(F^*+F\right)^{-\frac{1}{2}}
\left(F^*F+I_m\right)\left(F^*+F\right)^{-\frac{1}{2}}\|_2
\\~&=&
\sup\limits_{s\in\C_R}\rho\left(
\left(F^*+F\right)^{-\frac{1}{2}}
\left(F^*F+I_m\right)\left(F^*+F\right)^{-\frac{1}{2}}\right)
\\~&=&
\sup\limits_{s\in\C_R}\rho\left(
\left(F^*F+I_m\right)\left(F^*+F\right)^{-1}\right),
\end{matrix}
\]
so Eq. \eqref{eq:DefHPetaAgain2} is established.
\qed
\vskip 0.2cm

\noindent
The following alternative description of the set
$\mathcal{HP}_{\eta}$ will be employed in the sequel.

\begin{Cy}\label{Cy:AlternativeHPeta}
Let ${\scriptstyle\eta}$, ${\scriptstyle\eta}\in(1,~\infty]$,
be a parameter. An $m\times m$-valued rational function $F(s)$
belongs to $\mathcal{HP}_{\eta}$, if and only if, $F(s)$
analytically map $\C_R$ to $\mathbf{L}_{I_m}({\scriptstyle\eta})$,
see Eq. \eqref{eq:DefLh(eta)}.
\end{Cy}

\noindent
Indeed, from Observation \ref{Ob:CayleyQuantitativeHyperPositive1},
one can say that whenever $F(s)$ is an $m\times m$-valued
$\mathcal{HP}_{\eta}$ function, for each point $s_o\in\C_R$ the matrix
$F(s_o)$ belongs to $\mathbf{L}_{I_m}({\scriptstyle\eta})$. 
\vskip 0.2cm

\noindent
As an illustration, we first consider a structured, scalar,
$~\mathcal{HP}_{\eta}$ function of degree one.

\begin{Ex}\label{Ex:DegreeOne}
{\rm
In scalar terminology, Observation \ref{Ob:CayleyQuantitativeHyperPositive1}
may be interpreted as saying that that an $~\mathcal{HP}_{\eta}$ function,
maps $\C_R$ {\em to} a closed disk of the form,
\begin{equation}\label{eq:InvertibleDisk}
\mathbb{D}(\underbrace{\scriptstyle\eta}_{\rm Center},
\underbrace{\scriptstyle\sqrt{{\eta}^2-1}}_{\rm Radius}),
\quad
\quad
\quad
{\scriptstyle\eta}\in(1,~\infty].
\end{equation}
As a side remark, we mention that under inversion, this disk is mapped onto
itself. For further details see \cite[Section 3]{Lewk2020b}
\vskip 0.2cm

\noindent
Of particular interest, among the $~\mathcal{HP}_{\eta}$ functions, is
the degree one rational function of the form,
\begin{equation}\label{eq:ExampleStateSpaceDeg1diskPagain}
f(s)={\scriptstyle\eta}-{\scriptstyle\sqrt{\eta^2-1}}
+\frac{{\scriptstyle 2a}
{\scriptstyle\sqrt{\scriptstyle{\eta}^2-1}}}{s+{\scriptstyle a}}
\quad\quad\quad
\begin{smallmatrix}
\eta\in(1,~\infty]\\~\\
a>0,
\end{smallmatrix}
\end{equation}
which passes through the points
\mbox{$f(0)={\scriptstyle\eta}+{\scriptstyle\sqrt{{\eta}^2-1}}$},
\mbox{$~f(ia)={\scriptstyle\eta}-i{\scriptstyle\sqrt{{\eta}^2-1}}$}, 
and
\mbox{$~f(\infty)={\scriptstyle\eta}-{\scriptstyle\sqrt{{\eta}^2-1}}$}.
More generally, this $f(s)$ analytically maps $~\C_R~$ {\em onto} the
disk in Eq. \eqref{eq:InvertibleDisk}, irrespective of the value of the
parameter $a$. See the right-hand side of
Figure \ref{Figure:DiskCISs2}.
\vskip 0.2cm

\noindent
Recall that in Eq. \eqref{eq:DefScalarHyperBR} we considered the degree
one $\mathcal{HB}_{\eta}$ function
$g(s)={\scriptstyle\sqrt{\frac{\eta-1}{\eta+1}}}\frac{s-a}{s+a}$,
see the left-hand side of Figure \ref{Figure:DiskCISs2}. Recall
that in fact,
\[
f(s)=\mathcal{C}\left(g(s)\right)\quad\quad{\rm with}\quad\quad
\begin{smallmatrix}
g(s)~{\rm in~Eq.~}\eqref{eq:DefScalarHyperBR}
\\~\\
f(s)~{\rm in~Eq.~}\eqref{eq:ExampleStateSpaceDeg1diskPagain}.
\end{smallmatrix}
\]
As already mentioned following the work of O. Brune, the rational function
$f(s)$ in Eq. \eqref{eq:ExampleStateSpaceDeg1diskPagain} may be realized as
the driving point impedance of the simple electrical circuit in Figure
\ref{RLCofInvDisk}.

\begin{figure}[ht!]
\begin{tikzpicture}[scale=1.5]
  \draw[color=black, thick]
     (0,0) to [short,o-] (3,0){} % Baseline for connection to ground
     (-0.1,1) node[]{\large{$\mathbf{Z_{\rm in}~~\rightarrow}$}}
     (0,2) to [short,o-] (0.1,2)
     (0.1,2)  to [R,l=$\sqrt{\left(\frac{R}{2}\right)^2+1}-\frac{R}{2}$,](2,2)
     (2,2)   to node[short]{} (3,2)
     (3,0) to [C,l=$C$,] (3,2)
     (2,0) to [R, l=$R$, *-*] (2,2)
     ;
\end{tikzpicture}
\caption{$Z_{\rm in}={\scriptstyle\eta}-\sqrt{\scriptstyle{\eta}^2-1}+
\frac{2a\sqrt{{\eta}^2-1}}{s+a}\quad\quad C=\frac{1}{aR}~$.}
\label{RLCofInvDisk}
\end{figure}
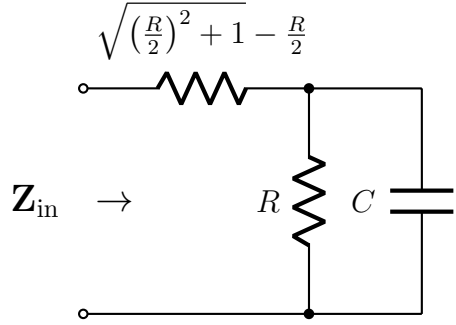

In the sequel, we shall relate to $f(s)$ in Eq.
\eqref{eq:ExampleStateSpaceDeg1diskPagain} as a representative of
$\mathcal{HP}_{\eta}$ functions of degree one.
}
\qed
\end{Ex}
\vskip 0.2cm

\subsection{The Lurie Problem - absolute stability}
\label{SubSec:Lurie}

\noindent
We have shown that within $\C_R$ the novel function set $\mathcal{HP}_{\eta}$ is
disk-contractive. Below we present an application of this observation to the
classical {\em Lurie problem} (a.k.a. the {\em absolute stability} problem). For
simplicity of exposition we here present the scalar version. For more information,
see e.g. \cite[Sections 3.2-3.5]{BroLozaMasEge2007},
\cite[Section 7.1]{Khalil2000}, \cite{Popov1973}.
\vskip 0.2cm

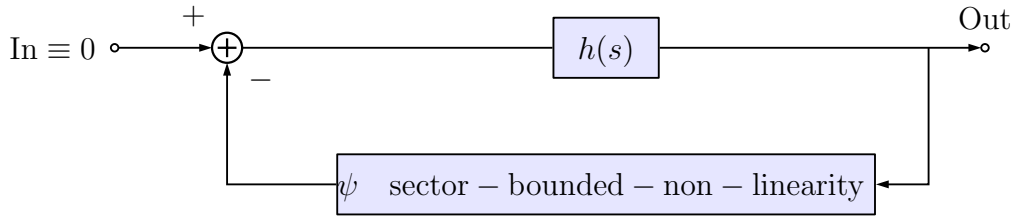
\begin{figure}
	 \begin{tikzpicture}
\matrix (diag_mat) [nodes={draw, fill=blue!10}, row sep=7mm, column sep=7mm]
		{
		\node[dspnodeopen,dsp/label=left] (x) {${\rm In}\equiv 0$};    &
	\node[dspadder,label={above left:$+$}, label={below right:$-$}] (sum) {}; &
		\node[dspfilter]    (G) {$~h(s)~$};&
		ֿ\coordinate (splt) {}; & 
		\node[dspnodeopen,dsp/label=above] (y) {${\rm Out}$};  \\
		~;& 		
		\coordinate (zr) ; &
		\node[dspfilter]    (H) {$\psi~~~{\rm sector-bounded-non-linearity}~$};&
		\coordinate (zl) {}; &
		~ ;\\ 
	};	
	\draw[dspconn] (x) -- (sum);
	\draw[dspconn] (sum) -- (G) -- (splt) -- (y);
	\draw[dspconn] (splt) -- (zl) -- (H);	
	\draw[dspconn] (H) -- (zr) -- (sum);	
	\end{tikzpicture}
\caption{Lurie problem feedback loop}
\label{NonLinearFeedbackLoop}
\end{figure}
\vskip 0.2cm

\noindent
{\bf Absolute stability problem:}~ Given a 
a feedback loop as in Figure \ref{NonLinearFeedbackLoop} where 
$~h(s)~$ is a rational function 
with neither pole nor zero at infinity and in the feedback,
\[
\psi=\psi(t, {\rm Out})
\]
is an {\em unknown time-dependent, sector-bounded, non-linearity},
satisfying
\begin{equation}\label{eq:SectorNonLinearity}
\left(K\cdot{\rm Out}-\psi\right)\left(\psi-k\cdot{\rm Out}\right)\geq 0
\quad\quad\quad
\begin{smallmatrix}
K,~k\in\R\quad K\geq k\\~\\
\forall~{\rm Out}\\~\\
\forall t\geq 0.
\end{smallmatrix}
\end{equation}
Does there exist $~K\geq k~$ so that the origin of the overall system is
uniformly asymptotically stable for any non-linearity\begin{footnote}{For
$~K=k$, this is a question on stability of~ {\em linear-time-invariant}
system, so we actually focus ourselves on the case $~K>k$.}\end{footnote}
$~\psi~$ satisfying Eq. \eqref{eq:SectorNonLinearity}.\\
Find sufficient (and possibly necessary) conditions for that and then
obtain various pairs $~K,~k~$ satisfying this requirement.
\qed
\vskip 0.2cm

\noindent
There are several absolute stability conditions see e.g.
\cite[Theorem 5.6.3]{AnderVongpa1973} \cite[Sections 7.1]{Khalil2000},
\cite{Popov1973}, \cite[Subsection 2.3.5]{SepJanKok1996}. 
We here refer only to the so called ``circle criterion"\begin{footnote}{
The graphical interpretation leading to the name ``circle criterion",
is beyond the scope of this work.}\end{footnote}.
\vskip 0.2cm

\noindent
{\bf The Circle Stability Criterion}\\

Consider the closed-loop system in Figure \ref{NonLinearFeedbackLoop}
with Eq. \eqref{eq:SectorNonLinearity}.\\
(I)
This system is absolutely stable whenever,
\begin{equation}\label{eq:LftCircleCriterion}
\left(f(h)\right)\in\mathcal{SP},
\end{equation}
for arbitrary rational function $f(s)$, of the form,
\begin{equation}\label{eq:DefGcircleCriterion}
f(s)={\scriptstyle\frac{1+Ks}{1+ks}}\quad{\rm with}\quad
K\geq k.
\end{equation}
(II) Whenever
$\infty >K\geq k>0$ the above  condition
can be equivalently formulated as:
This system is absolutely stable whenever,
\[
\left(f(h)\right)\in\mathcal{SP},
\]
for arbitrary rational function $f(s)$, of the form,
\begin{equation}\label{eq:DegOneHPetaAgain}
f(s)=
{\scriptstyle\left(\eta-\sqrt{{\eta}^2-1}\right)}+
{\scriptstyle\frac
{2a{\scriptstyle\sqrt{{\eta}^2-1}}}
{s+a}}
\quad{\rm with}\quad
\begin{matrix}
{\scriptstyle\eta}&=&{\scriptstyle\frac{1}{2}}\left(
\sqrt{\scriptstyle\frac{K}{k}}+\sqrt{\scriptstyle\frac{k}{K}}\right)
\\~\\
a&=&\frac{1}{K}~.
\end{matrix}
\end{equation}
\qed
\vskip 0.2cm

\noindent
Indeed, the statement in (I) is classical, see e.g.
\cite[Theorem 5.6.2]{AnderVongpa1973}, \cite[Theorem 7.1]{Khalil2000},
\cite[Sections 2.11, 3.10]{Rosenbr1974}, \cite[Proposition 3.6]{SepJanKok1996}. 
\vskip 0.2cm

\noindent
As for (II), recall that by item (ii) of Observation \ref{Ob:PropertiesP},
the set $\mathcal{SP}$ is a cone, and closed under inversion. Hence,
\[
f\left(h(s)\right)\in\mathcal{SP}
\quad\Longrightarrow\quad
\beta{f\left(h(s)\right)}^{-1}\in\mathcal{SP}
\quad\forall\beta>0.
\]
Now, restricting the discussion to $\infty >K\geq k>0$ 
means that one can take
\[
f(s)={\scriptstyle\beta}\left({\scriptstyle\frac{1+Ks}{1+ks}}\right)^{-1}=
{\scriptstyle\beta\frac{1+ks}{1+Ks}}_{|_{\beta=\frac{\sqrt{K}}{\sqrt{k}}}}
={\scriptstyle\frac{\sqrt{k}}{\sqrt{K}}+\frac{\left(
{\scriptstyle\frac{\sqrt{K}}{\sqrt{k}}}
-
{\scriptstyle\frac{\sqrt{k}}{\sqrt{K}}}
\right)\overbrace{\scriptstyle\frac{1}{K}}^a}
{s+\underbrace{\scriptstyle\frac{1}{K}}_{a}}}~,
\]
which is of the form of Eqs.  \eqref{eq:ExampleStateSpaceDeg1diskPagain}
and \eqref{eq:DegOneHPetaAgain}.
\vskip 0.2cm

\noindent
To summarize, the Circle Stability criterion may be formulated in
terms of a degree one rational function in $\mathcal{HP}_{\eta}$ as
in \eqref{eq:DegOneHPetaAgain}.

\subsection{Structure of $\mathcal{HP}_{\eta}$}

Following Eqs. \eqref{eq:LftCircleCriterion} and \eqref{eq:DegOneHPetaAgain},
we start by showing that the set $\mathcal{HP}_{\eta}$ is contractive in
the sense of~ {\em composition}.
\vskip 0.2cm

\noindent
In \cite{Reza1984} F.M. Reza, in the framework of scalar positive real
rational functions, coined the term ``power dominance": A function
$f_1(s)$ ``dominates" $f_2(s)$ whenever the Nyquist plot of $f_2(s)$
is contained in the Nyquist plot $f_1(s)$. Furthermore, under certain
conditions one can find a third function $h(s)$ so that these
functions are related through composition,
\[
f_2(s)=f_1\left(h(s)\right).
\]
Reza, subsequently applied this idea to the construction of $~R-L-C$ circuits.
\vskip 0.2cm

\noindent
\begin{Ex}\label{Ex:Reza}
{\rm
We here illustrate F.M. Reza's ``power dominance" in the framework
of $\mathcal{HP}_{\eta}$ functions. In particular, as the function $f(s)$
in Eq. \eqref{eq:DegOneHPetaAgain},
\[
f(s)=
{\scriptstyle\left(\eta-\sqrt{{\eta}^2-1}\right)}+
{\scriptstyle\frac
{2a{\scriptstyle\sqrt{{\eta}^2-1}}}
{s+a}}\quad\quad{\scriptstyle a}>0,
\]
analytically maps $~\C_R~$ {\em onto}
$~\mathbb{D}(\underbrace{\scriptstyle\eta}_{\rm Center},
\underbrace{\scriptstyle\sqrt{{\eta}^2-1}}_{\rm Radius})$,
this $f(s)$, dominates all functions in $\mathcal{HP}_{\eta}~$.
See the right-hand side of Figure \ref{Figure:DiskCISs2}
and for $\eta=\sqrt{2}$, the red curve in Figure \ref{Figure:Reza}.
\vskip 0.2cm

\noindent
(1)~
As was already pointed out, a degree one rational $\mathcal{HP}$
function is of the form
\begin{equation}\label{eq:ExDegOneHP}
\phi(s)=d+\frac{b}{s+a}\quad\quad\quad abd>0.
\end{equation}
This $\phi(s)$ analytically maps $~\C_R~$ to a closed disk
$~\mathbb{D}(
\underbrace{\scriptstyle d+\frac{b}{2a}+0i}_{\rm Center}~,
\underbrace{\scriptstyle\frac{b}{2a}}_{\rm Radius})~$
Substitution in Eq. \eqref{eq:DefHPetaAgain2} reveals that the function
in Eq. \eqref{eq:ExDegOneHP} in fact belongs to $\mathcal{HP}_{\eta}$
where,
\[
\eta=\left\{\begin{smallmatrix}\frac{1}{2}\left(d+\frac{1}{d}\right)&&&
d\in\left(0,~\frac{\sqrt{b^2+4a^2}-b}{2a}\right)
\\~\\
\frac{1}{2}\left(\frac{da+b}{a}+\frac{a}{da+b}\right)&&&
d>\frac{\sqrt{b^2+4a^2}-b}{2a}~.
\end{smallmatrix}\right.
\]
As already mentioned, all functions of the form of $\phi(s)$ in Eq.
\eqref{eq:ExDegOneHP} are dominated by the above $f(s)$.
\vskip 0.2cm

\noindent
(2)~ Consider now a degree two $\mathcal{HP}$ rational function
\begin{equation}\label{eq:DominatedFunction}
\phi(s)=\frac{\scriptstyle\sqrt{2}a^2}{(s+{\scriptstyle a})^2}+
{\scriptstyle\frac{1}{\sqrt{2}}}
\quad\quad\quad a>0.
\end{equation}
Note that $\phi(0)=\frac{3}{\sqrt{2}}$, $\phi(\pm ia)=\frac{1}{\sqrt{2}}(1\mp i)$,
$\phi(\pm{\scriptstyle\sqrt{3}}ai)=\frac{3}{4\sqrt{2}}\left(1\mp{i}\right)~$
and $\phi(\infty)=\frac{1}{\sqrt{2}}~$. In fact,
\[
\sup\limits_{s\in\C_R} 
\frac{{\scriptstyle\phi}^*{\scriptstyle\phi}+{\scriptstyle 1}}{{\scriptstyle\phi}^*+
{\scriptstyle\phi}}
=
{\frac{{\scriptstyle\phi}^*{\scriptstyle\phi}+{\scriptstyle 1}}{{\scriptstyle\phi}^*+
{\scriptstyle\phi}}}_{|_{s=\pm{i}a}}=\frac{2}{\frac{2}{\sqrt{2}}}=
{\scriptstyle\sqrt{2}},
\]
so from Eq. \eqref{eq:DefHPetaAgain2} one can conclude that $\phi(s)$ belongs to
$\mathcal{HP}_{\sqrt{2}}$, see the blue curve in Figure \ref{Figure:Reza}. As
already pointed out in the beginning of this
example, this $\phi(s)$ is dominated by,
\[
f(s)={\scriptstyle\sqrt{2}}-{\scriptstyle 1}+{\scriptstyle\frac{2a}{s+1}}~
\quad\quad{\scriptstyle a}>0,
\]
see the red curve in Figure \ref{Figure:Reza}. In particular note that
$
f\left(s\right)_{|_{s=\pm ia(\sqrt{2}+1)}}=
{\scriptstyle\frac{1}{\sqrt{2}}}(1\mp i)=\phi(s)_{|_{s=\pm ia}}.
$

\begin{figure}[ht!]
\includegraphics[width=0.70\textwidth]{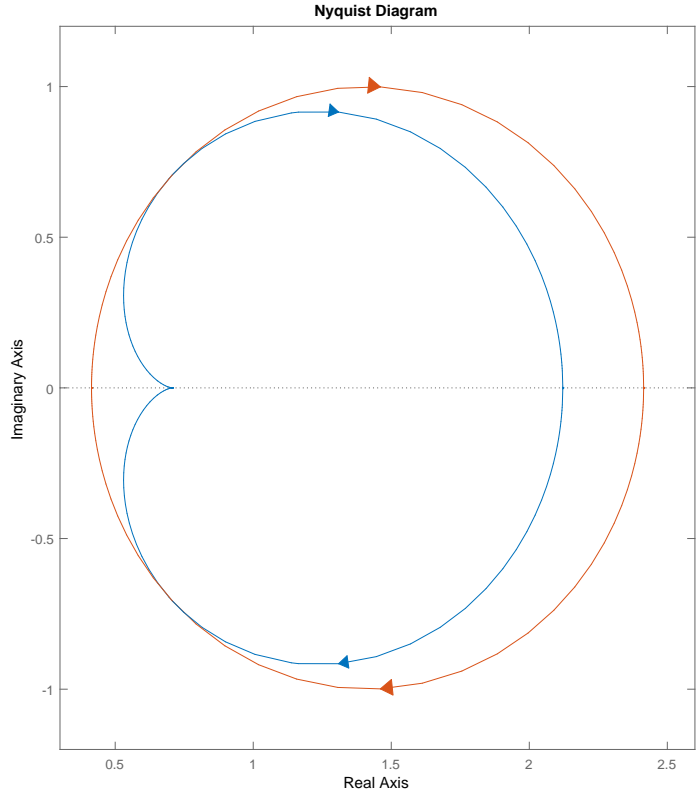}
\caption{Nyquist plots of $f(s)$ (red) and of $\phi(s)$ (blue),
from Eqs. \eqref{eq:ExampleStateSpaceDeg1diskPagain}
\eqref{eq:DominatedFunction}, respectively
}
\label{Figure:Reza}
\end{figure}
}
\qed
\end{Ex}
\vskip 0.2cm

\noindent
Next, we extend Reza's approach to a \mbox{matrix-valued} setup.

\begin{Ob}\label{Ob:Composition}
Let $F(s)$ be a \mbox{$m\times m$-valued} rational function in $\mathcal{HP}_{\eta}$
for some \mbox{${\scriptstyle\eta}\in(1, \infty]$}.  Let $h(s)$ be a scalar
rational function in $\mathcal{SP}$.\\
Then the composed \mbox{$m\times m$-valued} rational function
\[
F\left(h(s)\right)
\]
belongs to $\mathcal{HP}_{{\eta}_1}$ for some
\[
{\eta}_1\in(1, \eta].
\]
\end{Ob}
\vskip 0.2cm

\noindent
Indeed, as by definition $h(s)$ maps $\C_R$ {\em into} $\C_R$,
one has that,
\[
\begin{matrix}
{\scriptstyle\eta}=
\sup\limits_{s\in\C_R}\rho\left({\scriptstyle
\left(F(s)^*F(s)+I_m\right)\left(F(s)^*+F(s)\right)^{-1}}\right)
&\geq&
\sup\limits_{s\in h(\C_R)}\rho\left({\scriptstyle
\left(F(s)^*F(s)+I_m\right)\left(F(s)^*+F(s)\right)^{-1}}\right)
\\~&=&
\sup\limits_{s\in\C_R}\rho\left({\scriptstyle
\left(F(h)^*F(h)+I_m\right)\left(F(h)^*+F(h)\right)^{-1}}\right)
=
{\scriptstyle{\eta}_1}~.
\end{matrix}
\]
We further explore the structure of the set $\mathcal{HP}_{\eta}$.
For perspective recall that in \cite{Lewk2020a} it was shown that the
set of positive real functions is matrix-convex {\em cone}, closed
under inversion.

\begin{Pn}\label{Pn:StructureHPeta}
For a given ${\scriptstyle\eta}\in(1, \infty]$, the set  $\mathcal{HP}_{\eta}$ 
is matrix-convex. It is also closed under inversion, i.e. whenever a
rational function $F(s)$ is in $\mathcal{HP}_{\eta}$, then
$\left(F(s)\right)^{-1}$ is well defined and belongs
to $\mathcal{HP}_{\eta}$, with the same ${\scriptstyle\eta}$.
\end{Pn}

\noindent
{\bf Proof :}\quad Matrix-convexity. ~Recall that by item (II) of Theorem
\ref{Tm:ClosureOfL_H(r)} the set $\mathbf{L}_{I_m}({\scriptstyle\eta}$ is
matrix-convex. Combining this fact along with with the alternative
description of the set $\mathcal{HP}_{\eta}$ in Corollary
\ref{Cy:AlternativeHPeta}, establishes the claim.
\vskip 0.2cm

\noindent
We shall show invertibility in two ways.
\vskip 0.2cm

\noindent
(i) Direct approach. Using Eq. \eqref{eq:DefHPetaAgain2} one has that,
\[
\begin{matrix}
{\scriptstyle\eta}(F)&=&
\sup\limits_{s\in\C_R}\rho\left(
\left(F^*F+I_m\right)\left(F^*+F\right)^{-1}\right)
&~\\~&=&
\sup\limits_{s\in\C_R}
\rho\left(\underbrace{F^*\left(\left(F^*\right)^{-1}F^{-1}+I_m\right)
F}_{
\left(F^*F+I_m\right)
}
\underbrace{
F^{-1}\left(\left(F^*\right)^{-1}+F^{-1}\right)^{-1}\left(F^*\right)^{-1}
}_{\left(F^*+F\right)^{-1}}
\right)
&~\\~&=&
\sup\limits_{s\in\C_R}\rho
\left(
F^*
\left(\left(F^*\right)^{-1}F^{-1}+I_m\right)
\left(\left(F^*\right)^{-1}+F^{-1}\right)^{-1}(F^*)^{-1}
\right)&~\\~&=&
\sup\limits_{s\in\C_R}\rho
\left(
\left(\left(F^{-1}\right)^*F^{-1}+I_m\right)
\left(\left(F^{-1}\right)^*+F^{-1}\right)^{-1}
\right)&=
\eta(F^{-1}).
\end{matrix}
\]
(ii) Indirect (yet easy) approach. Applying Definition
\ref{Dn:MatrixCayleyTransform}, of the Cayley transform,
to a function $F(s)$, it is easy to verify that
\begin{equation}\label{eq:CayleyInvFunction}
\mathcal{C}\left(\left(F\right)^{-1}\right)
=-\mathcal{C}\left(F\right).
\end{equation}
Thus clearly, for any $s\in\C$ where $F(s)$ is analytic,
it holds that
\[
\|\mathcal{C}\left(F\right)\|_2
=
\|-\mathcal{C}\left(F\right)\|_2
=
\|\mathcal{C}\left(F^{-1}\right)\|_2~.
\]
This implies that
\[
\sup\limits_{s\in\C_R}\|\mathcal{C}\left(F(s)\right)\|_2=
\sup\limits_{s\in\C_R}\|\mathcal{C}\left(F(s)^{-1}\right)\|_2~,
\]
Assuming that $F(s)$ belongs to $\mathcal{HP}_{\eta}$, 
by Observation \ref{Ob:CayleyQuantitativeHyperPositive1}
it follows that this is equivalent to,
\[
{\scriptstyle\frac{\eta-1}{\eta+1}}=
\sup\limits_{s\in\C_R}\|\mathcal{C}\left(F(s)\right)\|_2
=
\sup\limits_{s\in\C_R}\|\mathcal{C}\left(F(s)^{-1}\right)\|_2~,
\]
so by Eq. \eqref{eq:NormBound}, the claim is established.
\qed

\section{Kalman-Yakubovich-Popov Lemma for $~\mathcal{HP}_{\eta}$ functions}
\label{Sect:K-Y-P}
\setcounter{equation}{0}

\noindent
A state-space realization of a given \mbox{$m\times m$-valued} rational
function $F(s)$, with no pole at infinity, is given by
\begin{equation}\label{eq:Realization}
F(s)=C(sI_n-A)^{-1}B+D\quad\quad\quad
R_F={\footnotesize\left(\begin{array}{c|c}
A&B\\ \hline C&D\end{array}\right)}.
\end{equation}
The $(n+m)\times(n+m)$ array $R_F$ was introduced by H.H. Rosenbrock,
see \cite{Rosenbr1974}. In general, $n$ is larger or equal to the
McMillan degree of $F(s)$. Equality means that the realization is minimal.
\vskip 0.2cm

\noindent
In the first subsection, we take advantage of the dual nature
of $R_F$ as an {\em array}~ and as a ~{\em matrix}.
\vskip 0.2cm

\subsection{K-Y-P for Positive functions}

\begin{Tm}\label{Tm:KYP}
Let $F(s)$ be a $m\times m$-valued rational function with no pole
at infinity and let $R_F$ in \eqref{eq:Realization} be a corresponding
$(n+m)\times(n+m)$ state-space realization.
\vskip 0.2cm

\noindent
For some $H\in\mathbf{P}_n$ consider the
following Lyapunov equation.
\begin{equation}\label{eq:Plemma}
\left(\begin{smallmatrix}-H&0\\0&I_m\end{smallmatrix}\right)
\left(\begin{smallmatrix}A&B\\
C&D\end{smallmatrix}\right)
+
\left(\begin{smallmatrix}A&B\\
C&D\end{smallmatrix}\right)^*
\left(\begin{smallmatrix}-H&0\\0&I_m\end{smallmatrix}\right)
=Q.
\end{equation}

\begin{itemize}
\item[(i)~~~] If 
in Eq. \eqref{eq:Plemma} $Q\in\overline{\mathbf P}_{n+m}$
then $F(s)$ is a $\mathcal{P}$ function.
\vskip 0.2cm

\noindent
\item[(ii)~~] If in Eq. \eqref{eq:Plemma} 
$~Q=Q_1+\left(\begin{smallmatrix}\Delta&0\\0&0\end{smallmatrix}\right)$
with $Q_1\in\overline{\mathbf P}_{n+m}$ and
${\scriptstyle\Delta}\in\mathbf{P}_n$, then $F(s)$ is a $~\mathcal{SP}$ function.
\vskip 0.2cm

\noindent
\item[(iii)~] If in Eq. \eqref{eq:Plemma} $~Q\in{\mathbf P}_{n+m}$,
then $F(s)$ is a $\mathcal{HP}$ function.
\end{itemize}
\vskip 0.2cm

\noindent
In each of the three above statements, if the realization in Eq.
\eqref{eq:Realization} is minimal, i.e.  $n$ is the McMillan degree, then
the converse is true as well.
\end{Tm}

\noindent
{\bf Proof :}\quad (i) Up to a change of coordinates one can substitute in
Eq. \eqref{eq:Plemma} $H=I_n$ and thus,
\[
Q=
\left(\begin{smallmatrix}-I_n&0\\0&I_m\end{smallmatrix}\right)
\underbrace{
\left(\begin{smallmatrix}A&B\\
C&D\end{smallmatrix}\right)}_{R_F}+
\underbrace{
\left(\begin{smallmatrix}A&B\\
C&D\end{smallmatrix}\right)^*}_{R_F^*}
\left(\begin{smallmatrix}-I_n&0\\0&I_m\end{smallmatrix}\right)
=
\left(\begin{smallmatrix}-A&-B\\
~~C&~~D\end{smallmatrix}\right)
\left(\begin{smallmatrix}-A&-B\\
~~C&~~D\end{smallmatrix}\right)^*
\]
Next note that for all $s\in\C$,
\[
\begin{matrix}
\hat{Q}(s):=
\left(\begin{smallmatrix}2{\rm Re}(s)I_n&0\\
0&0\end{smallmatrix}\right)
+
Q
&=&
\left(\begin{smallmatrix}sI_n&0\\
0&0\end{smallmatrix}\right)
+
\left(\begin{smallmatrix}-A&-B\\
~~C&~~D\end{smallmatrix}\right)
+
\left(\begin{smallmatrix}sI_n&0\\
0&0\end{smallmatrix}\right)^*
+
\left(\begin{smallmatrix}-A&-B\\
~~C&~~D\end{smallmatrix}\right)^*
\\~\\~&=&
\left(\begin{smallmatrix}sI_n-A&&-B\\
~~C&&~~D\end{smallmatrix}\right)
+
\left(\begin{smallmatrix}sI_n-A&&-B\\
~~C&&~~D\end{smallmatrix}\right)^*.
\end{matrix}
\]
Now,
\[
\left(\begin{smallmatrix}\left((C(A-sI_n)^{-1}\right)^*\\~
\\I_m\end{smallmatrix}\right)^*\hat{Q}(s)\left(\begin{smallmatrix}
\left((C(A-sI_n)^{-1}\right)^*\\~\\I_m\end{smallmatrix}\right)
=\underbrace{\scriptstyle C(sI_n-A)^{-1}B+D}_{F(s)~{\rm in~Eq.}~
\eqref{eq:Realization}}+\underbrace{\scriptstyle\left(C(sI_n-A)^{-1}B
+D\right)^*}_{\left(F(s)\right)^*~{\rm in~Eq.}~\eqref{eq:Realization}}
\]
and if $Q\in\overline{\mathbf P}_{n+m}$ (or
$Q\in\mathbf{P}_{n+m}$) then for all $s\in\C_R$
one has that $\hat{Q}(s)\in\overline{\mathbf P}_{n+m}$ (or
$\hat{Q}(s)\in\mathbf{P}_{n+m}$) and then
$F(s)\in\overline{\mathbf L}_{I_m}$
(or $F(s)\in\mathbf{L}_{I_m}$) respectively. Hence items (i) and (iii)
are established.
\vskip 0.2cm

\noindent
(ii) By definitions $F\in\mathcal{SP}$ means that for some $\epsilon>0$ one has that 
\mbox{$C\left(sI_n-(A+\epsilon{I}_n)\right)^{-1}B+D$} is positive real.
Namely, by item (i), for some $~H\in\mathbf{P}_n$,
\begin{equation}\label{eq:SPlemma}
\left(\begin{smallmatrix}-H&0\\0&I_m\end{smallmatrix}\right)
\left(\begin{smallmatrix}A+\epsilon{I}_n&B\\
C&D\end{smallmatrix}\right)+
\left(\begin{smallmatrix}A+\epsilon{I}_n&B\\
C&D\end{smallmatrix}\right)^*
\left(\begin{smallmatrix}-H&0\\0&I_m\end{smallmatrix}\right)
=Q\in\overline{\mathbf{P}}_{n+m}~.
\end{equation}
Then,
\[
\left(\begin{smallmatrix}-H&0\\0&I_m\end{smallmatrix}\right)
\left(\begin{smallmatrix}A&B\\
C&D\end{smallmatrix}\right)+
\left(\begin{smallmatrix}A&B\\
C&D\end{smallmatrix}\right)^*
\left(\begin{smallmatrix}-H&0\\0&I_m\end{smallmatrix}\right)
=\underbrace{Q}_{\in\overline{\mathbf{P}}_{n+m}}
+\left(\begin{smallmatrix}
2\epsilon{H}&0\\0&0\end{smallmatrix}\right),
\]
so taking
$\left({\scriptstyle 2\|H\|_2\|{\scriptstyle\Delta}^{-1}\|_2}\right)^
{-1}\geq\epsilon$
completes the construction.
\qed
\vskip 0.2cm

\noindent
There are numerous proofs for variants of Theorem \ref{Tm:KYP}, see e.g.
\cite{AlpayLew2011}, \cite[Chapter 5]{AnderVongpa1973},
\cite[Lemma 6.3, Appendix C.12]{Khalil2000}. The formulation in Eq.
\eqref{eq:Plemma} (introduced in \cite{DickDelsGenKam1985}) enables
one to obtain the above simple proof.
\vskip 0.2cm

\noindent
In \cite[Section 5]{Lewk2020a} we advance in manipulating $R_F$ in
Eq. \eqref{eq:Realization} as a matrix, and consider families of
rational functions whose corresponding realization arrays/matrices
$R_F$, form a convex cone, closed under inversion.
\vskip 0.2cm

\noindent
As a consequence of Theorem \ref{Tm:KYP}, we have the following.

\begin{Ob}\label{Ob:ThreeVersionsHyperPosKyp}
Let $F(s)$ be a $m\times m$-valued rational function. Then
$F\in\mathcal{HP}$ if and only if, $F\in\mathcal{SP}$ and in addition,
$\lim\limits_{s~\rightarrow~\infty}F(s)$ exists and belongs to
$\mathbf{L}_{I_m}$.
\end{Ob}

\noindent
{\bf Proof : :}\quad Following Theorem \ref{Tm:KYP}, we actually need
to show that for a given $H\in\mathbf{P}_n$ in Eq. \eqref{eq:Plemma},
having on the right-hand side $Q\in\mathbf{P}_{n+m}$ is equivalent to
\[
\left(\begin{smallmatrix}-H&0\\0&I_m\end{smallmatrix}\right)
\left(\begin{smallmatrix}A+\epsilon{I}_n&B\\
C&D-\delta{I}_m\end{smallmatrix}\right)+
\left(\begin{smallmatrix}A+\epsilon{I}_n&B\\
C&D-\delta{I}_m\end{smallmatrix}\right)^*
\left(\begin{smallmatrix}-H&0\\0&I_m\end{smallmatrix}\right)
=\hat{Q}\quad\quad\hat{Q}\in\overline{\mathbf P}_{n+m}~,
\]
for some $~\epsilon$, $\delta>0$. Indeed, comparison yields
\[
\underbrace{Q}_{\in\mathbf{P}_{n+m}}=
\underbrace{\hat{Q}}_{\in\overline{\mathbf P}_{n+m}}+
{\scriptstyle 2}\left(\begin{smallmatrix}
{\scriptstyle\epsilon}H&0\\
0&{\scriptstyle\delta}I_m
\end{smallmatrix}\right),
\]
so having $\hat{Q}\in\overline{\mathbf P}_{n+m}$ implies that
$Q\in\mathbf{P}_{n+m}~$.
\vskip 0.2cm

\noindent
Conversely, by assumption
\[
\beta:={\scriptstyle\frac{1}{2}}\min\limits_{j=1,~\ldots~,~n+m}\lambda_j(Q),
\]
is positive. Thus taking,
\[
\epsilon\leq\beta\|H\|_2^{-1}
\quad\quad{\rm and}\quad\quad
\delta\leq\beta,
\]
guaranties that indeed $\hat{Q}\in\overline{\mathbf P}_{n+m}~$.
\vskip 0.2cm

\noindent
The above reasoning guarantees that $\delta>0$ and thus indeed,
$D\in\mathbf{L}_{I_m}$, so the claim is established.
\qed
\vskip 0.2cm

\subsection{K-Y-P for Hyper-Positive functions}
\label{SubSec:K-Y-PforHyper-Pos}

Next we consider $m\times m$-valued rational function $F(s)$
assuming it has no pole at infinity, and is of
McMillan degree $n$. Thus, it admits a state space realization
\begin{equation}\label{eq:RealizationF4}
R_F=\left({\footnotesize\begin{array}{l|r}
A&B\\
\hline
C&D
\end{array}}\right).
\end{equation}
Recall that applying the Cayley transform (see Definition
\ref{Dn:MatrixCayleyTransform}) to $F(s)$ means that
whenever $-1\not\in{\rm spect.}(D)$,
\begin{equation}\label{eq:RealizationUnderCayleyTrans}
R_{{\mathcal C}(F)}=
R_{2(F+I_m)^{-1}-I_m}=
\left({\footnotesize\begin{array}{l|r}
\overbrace{A-B\left(I_m+D\right)^{-1}C}^{\hat{A}}
&
\overbrace{
-{\scriptstyle\sqrt{2}}B\left(I_m+D\right)^{-1}}^{\hat{B}}\\
\hline
\underbrace{
{\scriptstyle\sqrt{2}}\left(I_m+D\right)^{-1}C}_{\hat{C}}&
\underbrace{
\left(I_m+D\right)^{-1}\left(I_m-D\right)}_{\hat{D}}
\end{array}}\right)
\end{equation}
Substituting Eq.  \eqref{eq:RealizationUnderCayleyTrans} in Eq.
\eqref{eq:QuantitativeBRlemmaAlternative} yields
the main result of this subsection.

\begin{Tm}\label{Tm:KypEtaHyperPositive}
Let $F(s)$ be a $m\times m$-valued rational function with no
pole at infinity of McMillan degree $n$ whose minimal realization
is given in Eq. \eqref{eq:RealizationF4}.
\vskip 0.2cm

\noindent
For some ${\scriptstyle\eta}\in(1,~\infty]$, the function $F(s)$ belongs
to the class $\mathcal{HP}_{\eta}$, if and only if, there exists
a matrix $H$ satisfying,
\[
H\in\mathbf{P}_n~,
\]
and
\[
\left(\begin{smallmatrix}
A-B\left(I_m+D\right)^{-1}C& -{\scriptstyle\sqrt{2}}B\left(I_m+D\right)^{-1}\\
{\scriptstyle\sqrt{2}}\left(I_m+D\right)^{-1}C&(I_m+D)^{-1}(I_m-D)\\
I_n&0\\
0&I_m
\end{smallmatrix}\right)^*
\underbrace{
\left(\begin{smallmatrix}
~~0&0&-H&~0\\
~~0&{\scriptstyle\frac{1+\eta}{1-\eta}}I_m&~~0&~0\\
-H&0&~~0&~0\\
~~0&0&~~0&~I_m
\end{smallmatrix}\right)
}_W
\left(\begin{smallmatrix}
A-B\left(I_m+D\right)^{-1}C& -{\scriptstyle\sqrt{2}}B\left(I_m+D\right)^{-1}\\
{\scriptstyle\sqrt{2}}\left(I_m+D\right)^{-1}C&(I_m+D)^{-1}(I_m-D)\\
I_n&0\\
0&I_m
\end{smallmatrix}\right)
\in\overline{\mathbf P}_{n+m}~.
\]
\end{Tm}
\vskip 0.2cm

\noindent
As already mentioned, the Hermitian matrix $W$, has $n+m$ positive and $n+m$
negative eigenvalues.
\vskip 0.2cm

\noindent
We next illustrate the quantitative essence of Theorem
\ref{Tm:KypEtaHyperPositive}.
\vskip 0.2cm

\begin{Ex}\label{Ex:KypQuantitativeHyperStable}
{\rm
Substituting in Eq. \eqref{eq:ExampleStateSpaceDeg1diskPagain} $\eta=\frac{5}{3}$
and $a=\frac{1}{9}$, one obtains the degree one $\mathcal{HP}_{\frac{5}{3}}$ function
$~f(s)={\scriptstyle\frac{1}{3}}+\frac{\frac{8}{27}}{s+\frac{1}{9}}=
\frac{\frac{1}{3}\left(s+1\right)}{s+\frac{1}{9}}~$. Recall now that
\begin{equation}\label{eq:f1}
f_1(s)={\scriptstyle\left(\frac{1}{2}\left(f(s)+\frac{1}{f(s)}\right)\right)^{-1}}
=
{\scriptstyle\frac{3}{5}}+\frac{
{\scriptstyle\frac{32}{75}}s}{s^2+{\scriptstyle\frac{2}{5}}s+{\scriptstyle\frac{1}{9}}}
={\scriptstyle\frac{3(s^2+\frac{10}{9}s+\frac{1}{9})}
{5(s^2+\frac{2}{5}s+\frac{1}{9})}}~,
\end{equation}
is a degree two function in $\mathcal{HP}_{{\eta}_1}$ where
\[
{\eta}_1={\scriptstyle\frac{1}{2}}\left({\scriptstyle\eta}+
{\scriptstyle\frac{1}{\eta}}\right)=
{\scriptstyle\frac{1}{2}}\left({\scriptstyle\frac{5}{3}}
+{\scriptstyle\frac{3}{5}}\right)
={\scriptstyle\frac{17}{15}}~.
\]
In particular, $~f_1(0)={\scriptstyle\frac{3}{5}}$, 
$~f_1({\scriptstyle\pm\frac{i}{3}})={\scriptstyle\frac{5}{3}}$,
and
$~f_1(\infty)={\scriptstyle\frac{3}{5}}$.
\vskip 0.2cm

\noindent
As a side remark, note that applying the Cayley transform here
yields the corresponding $\mathcal{HB}_{\eta}$ functions,
\mbox{$\mathcal{C}\left(f(s)\right)=g(s)=
{\scriptstyle\frac{1}{2}}\cdot\frac{s-\frac{1}{3}}{s+\frac{1}{3}}$}
and
\mbox{$\mathcal{C}\left(f_1(s)\right)=\left(g(s)\right)^2=g_1(s)=
{\scriptstyle\frac{1}{4}}\frac{(s-\frac{1}{3})^2}{(s+\frac{1}{3})^2}~$}.
\noindent

\vskip 0.2cm
Now, a balanced realization array of $f_1(s)$ in Eq.  \eqref{eq:f1} is given by,
\begin{equation}\label{eq:ExampleBalancedRealizationDegTwoSPDP}
R_{f_1}={\scriptstyle\frac{1}{5}}{\footnotesize\left(\begin{array}{rr|r}
-1&-\frac{4}{3}&\frac{4}{\sqrt{6}}\\
\frac{4}{3}&-1&\frac{8}{\sqrt{6}}\\
\hline
\frac{8}{\sqrt{6}}&\frac{4}{\sqrt{6}}&3
\end{array}\right)}.
\end{equation}
As before, from the classical K-Y-P Lemma one can (only qualitatively)
conclude that $f_1(s)$ in Eq. \eqref{eq:f1}
is Hyper-Positive. Indeed, taking in Eq. \eqref{eq:Plemma},
\[
H=I_2~,
\]
yields the following right-hand side,
\[
Q={\scriptstyle\frac{2}{5}}\left(\begin{smallmatrix}
1&0&\frac{2}{\sqrt{6}}\\
0&1&-\frac{2}{\sqrt{6}}\\
\frac{2}{\sqrt{6}}&-\frac{2}{\sqrt{6}}&3
\end{smallmatrix}\right),
\]
which is positive definite.
\vskip 0.2cm

\noindent
To {\em quantitatively}~ examine this function, we resort to Theorem
\ref{Tm:KypEtaHyperPositive}. Following Eq.
\eqref{eq:RealizationUnderCayleyTrans}, a minimal realization
of $~{\mathcal C}(f_1)$, the image of this $f_1(s)$ under the
Cayley transform, is,
\[
R_{{\mathcal C}(f_1)}=
{\footnotesize\left(\begin{array}{rr|c}
-\frac{1}{3}&-\frac{1}{3}&-\frac{1}{2\sqrt{3}}\\
0&-\frac{1}{3}&-\frac{1}{\sqrt{3}}\\
\hline
\frac{1}{\sqrt{3}}&\frac{1}{2\sqrt{3}}&\frac{1}{4}
\end{array}\right)}.
\]
Now taking in Lemma \ref{La:QuantitativeBRlemma}
(or in Theorem \ref{Tm:KypEtaHyperPositive}) 
\[
{\scriptstyle\eta}_1
={\scriptstyle\frac{17}{15}}\quad{\rm and}\quad
H=\left(\begin{smallmatrix}8&0\\0&2\end{smallmatrix}\right),
\]
yields that the right-hand side of Eq. 
\eqref{eq:QuantitativeBRlemmaAlternative} (or of the quadratic inclusion
in Theorem \ref{Tm:KypEtaHyperPositive}) is $~0_{3\times 3}~$, so indeed
$f_1\in\mathcal{HP}_{\frac{17}{15}}~$.
}
\qed
\end{Ex}


\vskip 0.2cm


\begin{thebibliography}{ZZ}

%\bibitem{AlpGoh1998}D. Alpay and I. Gohberg,~ ``Unitary Rational Matrix
%Functions" In I. Gohberg, editor, Topics in Interpolation Theory of
%Rational Matrix-Valued Functions,~ {\em Operator Theory: Advances and
%Applications},~ Vol. 33, pp. 175-222. Birkh{\" a}user Verlag, Basel, 1988.
%
\bibitem{AlpayLew2011}D. Alpay and I. Lewkowicz,~ ``The Positive Real Lemma and
Construction of all Realizations of Generalized Positive Rational Functions",~
{\em Systems and Control Letters},~ Vol. 60, pp. 985-993, 2011.

%\bibitem{AlpayLew2013}D. Alpay and I. Lewkowicz, ~``Convex Cones of
%Generalized Positive Rational Functions and the Nevanlinna-Pick
%Interpolation", ~{\em Linear Algebra and its Applications},
%Vol. 438, pp. 3949-3966, 2013.
%
\bibitem{AlpayLew2019a}D. Alpay and I. Lewkowicz,~ ``Composition of Rational
Functions: State-space Realization and Applications",  to appear in
{\em Linear Algebra and its Applications}.  See~ arXiv:1807.01753.

\bibitem{AlpLew2019b}D. Alpay and I. Lewkowicz,~ ``Realization of
Tensor-Product and of Tensor-Factorization of Rational Functions",~
{\em Quantum Studies: Mathematics and Foundations},~ Vol. 6,
pp. 269-278, 2019.

\bibitem{AnderVongpa1973}B.D.O. Anderson and S. Vongpanitlerd,~ {\em Networks
Analysis and Synthesis, A Modern Systems Theory Approach,}~ Prentice-Hall,
New Jersey, 1973.

\bibitem{Ando2004}T. Ando,~ ``Sets of Matrices with Common Stein Solutions
and H-contractions”,~ {\em Linear Algebra and its Application},~ Vol. 383,
pp. 49-64, 2004.

\bibitem{Belev1968}V. Belevich,~ {\em Classical Network Theory},~
Holden Day, San-Francisco, 1968.

\bibitem{BGFB1994}S. Boyd, L. El-Ghaoui, E. Ferron and
V. Balakrishnan, ~{\em Linear Matrix Inequalities in
Systems and Control Theory}, SIAM books, 1994.

\bibitem{BroLozaMasEge2007}B. Brogliato, R. Lozano, B. Maschke
and O. Egeland,~ {\em Dissipative Systems Analysis and
Control: Theory and Applications},~ Second Edition, Upper Saddle River,
New-York, Springer Verlag, 2007.

\bibitem{Brune1}O. Brune, ``Synthesis of a Finite Two
Terminal Network whose Driving Point Impedance is a
Prescribed Function of Frequency",~ {\em Journal of
Mathematical Physics},~ Vol. 10, pp. 191-236, 1931.

\bibitem{CohenLew1997a}N. Cohen and I. Lewkowicz,~ ``Convex Invertible
Cones and the Lyapunov Equation", ~{\em Linear Algebra
and its Applications},~ Vol. 250, pp. 265-286, 1997.
 
%\bibitem{CohenLew2007}N. Cohen ~and~ I. Lewkowicz,~ ``Convex Invertible
%Cones and Positive Real Analytic Functions", ~{\em Linear Algebra
%and its Applications},~ Vol. 425, pp. 797-813, 2007.
%
\bibitem{DickDelsGenKam1985}B. Dickinson, Ph. Delsarte, Y. Genin and
Y. Kamp,~ ``Minimal Realization of Pseudo Positive and Pseudo
Bounded Real Rational Matrices”,~ {\em IEEE trans. Circuits and Systems},~
Vol. 32, pp. 603-605, 1985.

%\bibitem{HornJohnson1}R. A. Horn and C. R. Johnson,~ {\em Matrix ~
%Analysis},~ Cambridge University Press, 1985.
%
%\bibitem{HornJohnson2}R. A. Horn and C. R. Johnson,~ {\em Topics
%in Matrix ~Analysis},~ Cambridge University Press, 1991.
%
%\bibitem{Kailath1980}T. Kailath,~ {\em Linear Systems},~
%Prentice-Hall, 1980.
%
\bibitem{Khalil2000}H.K. Khalil, ~{\em Nonlinear Systems},~
$3^{\rm rd}$ edition, Pearson Education, NJ, USA, 2000.

\bibitem{Lewk2020a}I. Lewkowicz,~ {\em Passive Linear Systems Matrix-convex
Invertible Cones Point of View},~ a manuscript.

\bibitem{Lewk2020b}I. Lewkowicz,~ {\em Matrix Sign Function Iterations -
Geometric Point of View},~ a manuscript.

\bibitem{Lewk2020c}I. Lewkowicz,~ {\em Discrete-Time Passive Linear Systems -
Multiplicative Matrix-Convex Sets Point of View},~ a manuscript.

\bibitem{MorelliSmith2019}A. Morelli and M.C. Smith, {\em
Passive Network Synthesis: An Approach to Classification}
no. DC33 in {\em Advances in Design and Control} series
by SIAM, 2019.

\bibitem{Popov1973}V.M. Popov,~ {\em Hyper-stability of Control
Systems},~ Springer 1973.

%\bibitem{Rant1993}A. Rantzer, ~``A ~Weak ~Kharitonov ~Theorem ~Holds if
%and only if the Stability ~Region and its Reciprocal are Convex",~{
%\it International Journal of Robust and Nonlinear Control},~ Vol. 3,
%pp. 55-62, 1993.
%
\bibitem{Reza1984}F.M. Reza,~ ``The Concept of Power Dominant
Systems",~ {\em Lecture Notes in Control and Information
Sciences} Vol. 58, pp. 787-795, Springer 1984.

\bibitem{Rosenbr1974}H.H. Rosenbrock,~ {\em Computer-Aided ~Control~
System Design},~ Academic ~Press, 1974.

\bibitem{SepJanKok1996}R. Sepulchre, M. Jankovi\'{c} and P.V.
Kokotovi\'{c},~{\em Constructive Nonlinear Control},~
Communication and Control Engineering series, Springer, 1996.

\bibitem{Will1976}J.C. Willems,~ ``Realization of Systems with Internal
Passivity and Symmetry Constraints",~ {Journal of the Franklin
Institute}, Vol. 301, pp. 605-621, 1976.

\bibitem{Wohl1969}M. R. Wohlers,~ {\em Lumped and Distributed
Passive Networks},~ Acad. Press 1969.

\end{thebibliography}
\end{document}